\documentclass[10pt,reqno]{article}
\usepackage[top=1in, bottom=1in, left=1.2in, right=1.2in, asymmetric]{geometry}

\usepackage{fancyhdr}
\usepackage{amsmath,amsthm,amssymb,fancybox,ifthen,float}
\usepackage{graphicx}
\usepackage{setspace}
\usepackage{comment}
\usepackage[T1]{fontenc}
\usepackage{mathpazo}
\usepackage{sectsty}
\usepackage{bm}

\usepackage{subcaption}
\usepackage{color}

\sectionfont{\large}
\subsectionfont{\normalsize}
\numberwithin{equation}{section}

\theoremstyle{definition}

\theoremstyle{remark}

\floatstyle{plain}
\restylefloat{figure}

\newcommand\myout{\text{out}}
\newcommand\myin{\text{in}}

\newcommand{\lth}{l^{\text{th}}}

\newcommand\Chi{\mathcal{X}}

\newcommand\cK{\mathcal K}

\newcommand\bc{\boldsymbol c}

\newcommand\bphi{\boldsymbol \phi}

\newcommand\btau{\boldsymbol \tau}

\newcommand\bj{\boldsymbol j}

\newcommand\bmm{\boldsymbol m}

\newcommand\bu{\boldsymbol u}

\newcommand\bx{\boldsymbol x}

\newcommand\bA{\boldsymbol A}
\newcommand\bQ{\boldsymbol Q}
\newcommand\ba{\boldsymbol a}
\newcommand\bB{\boldsymbol B}

\newcommand\bH{\boldsymbol H}

\newcommand\bl{\boldsymbol l}
\newcommand\bE{\boldsymbol E}
\newcommand\bF{\boldsymbol F}
\newcommand\bJ{\boldsymbol J}

\newcommand\cS{\mathcal{S}}

\renewcommand\Re{\operatorname{Re}}

\newcommand{\surfgrad}{\nabla_\Gamma }
\newcommand{\surflap}{\Delta_\Gamma}

\newcommand{\vct}[1]{\bm{\mathsf{#1}}}

\newcommand{\mtx}[1]{\bm{\mathsf{#1}}}

\newcommand{\fluxtor}{\Phi^{\text{tor}}}
\newcommand{\fluxpol}{\Phi^{\text{pol}}}

\makeatletter
\let\@fnsymbol\@alph
\makeatother

\begin{document}

\title{An integral equation-based numerical solver for 
Taylor states in toroidal geometries}

\author{Michael O'Neil\footnote{Courant Institute and Tandon
School of Engineering, New York
    University, New York.  E-mail:
    {\tt oneil@cims.nyu.edu}. Research of M.~O'Neil partially supported by
    the Air Force Office of Scientific Research under NSSEFF Program
    Award FA9550-10-1-0180 and the Office of Naval Research under
    awards N00014-17-1-2059 and N00014-17-1-2451.
    (Corresponding author.)} \, and
  Antoine J. Cerfon\footnote{Courant
    Institute, New York University, 251 Mercer Street, New York, NY.  
    E-mail: {\tt cerfon@cims.nyu.edu}. Research of A.~J.~Cerfon
    partially supported by U.S. Department of Energy, Office of Science, Fusion Energy Sciences under Awards No. DE-FG02-86ER53223
and  DE-SC0012398}}

\maketitle

\begin{abstract}
  We present an algorithm for the numerical calculation of Taylor
  states  in toroidal
  and toroidal-shell geometries using an analytical framework developed
  for the solution to the time-harmonic Maxwell equations. Taylor
  states are a special case of what are known as Beltrami fields,
  or linear force-free fields. The scheme of this work
  relies on the generalized Debye source representation of Maxwell
  fields and an integral representation of Beltrami fields which immediately
  yields a well-conditioned second-kind integral equation. This
  integral equation has a unique solution whenever the Beltrami
  parameter $\lambda$ is not a member of a discrete, countable set of
  resonances which physically correspond to spontaneous symmetry
  breaking.  Several numerical examples relevant to
  magnetohydrodynamic equilibria calculations are
  provided. Lastly, our approach easily generalizes to arbitrary geometries,
  both bounded and unbounded, and of varying genus.\\

  \noindent {\bf Keywords}: Beltrami field, Beltrami flow, generalized
  Debye sources, time-harmonic Maxwell's equations, plasma physics,
  electromagnetics, Debye sources, force-free fields, Taylor states,
  magnetohydrodynamics
\end{abstract}

\onehalfspacing

\section{Introduction}

A wide range of astrophysical and laboratory plasmas are in force-free
equilibria \cite{flyer,amari1,ram,taylor_review} where the
magnetic field $\bB$ satisfies
\begin{equation}\label{eq-forcefree}
(\nabla \times \bB) \times \bB = \boldsymbol{0}
\end{equation}
This immediately implies that the current density
$\bJ=(\nabla \times \bB)/\mu_{0}$ is parallel to  the magnetic
field, i.e. there exists a scalar function
$\lambda = \lambda(\boldsymbol{x})$ such that
\begin{equation}\label{eq-linearity}
\nabla\times \bB =\lambda  \bB.
\end{equation}
Within the general class of linear force-free equilibria described
by~\eqref{eq-linearity}, Taylor states or Woltjer-Taylor states are
particular equilibria for which $\lambda$ is a spatially uniform
constant given by the ratio of the magnetic energy to the magnetic
helicity, see Chapter 11 of~\cite{bellantext}.
They play a central role in plasma
physics~\cite{chandrasekhar,woltjer,bruno,tang,cothranprl,
  battaglia,hudson,dennis}
as the natural state resulting from dissipative turbulent
relaxation~\cite{taylor_review,qin}. Since they satisfy the equation
$\nabla\times \bB =\lambda \bB$ with $\lambda$ constant, magnetic
fields in Taylor state configurations are a special case of a class of
force-free fields called linear Beltrami
fields~\cite{amari2,boulmezaoud}. Since in this work we will consider
$\lambda$ to be a given input to the solver, we will use the
expressions \emph{Taylor state} and \emph{linear Beltrami field}
interchangeably for the remainder of this article.

Linear Beltrami fields have been extensively studied mathematically, and their properties are well
understood by now~\cite{hiptmair,picard,EpGrOn2,kress_1986,constantin_1988}. On the
other hand, relatively few numerical solvers have been developed to
compute them in geometries relevant to plasma physics. 
To the best of the authors' knowledge, solvers for these problems
based on integral equation formulations have never been constructed
despite their desirable properties: access to relatively
high-precision derivatives of the field (via analytic differentiation
of the integral representation), low memory requirements (only the
boundary has to be discretized), and overall rapid convergence of the
solution (when coupled with high-order quadrature rules and a fast
algorithm, such as a fast multipole method).

We take a moment to justify this claim. One may at first think that
Taylor states in axisymmetric toroidal geometries can be viewed as a
special class of more general Grad-Shafranov
equilibria~\cite{grad,shafranov}, as was for instance done in
\cite{cerfontaylor}. From this point of view, Grad-Shafranov solvers
relying on integral formulations~\cite{lao,pataki} could be used to
compute linear Beltrami fields. However, this approach is not
satisfactory for the following reasons. First, a Grad-Shafranov
solver would not take advantage of the particular properties of
linear Beltrami fields. Second, certain applications~\cite{bruno,hudson}
require the computation of linear Beltrami fields in hollow toroidal
shells. Grad-Shafranov solvers are usually not designed to handle such
geometries. Finally, and most importantly, Taylor states in
axisymmetric domains may not be axisymmetric themselves
\cite{taylor_review}. By definition, the Grad-Shafranov equation does
not apply to these fully three-dimensional, bifurcated states.

The purpose of this article is to present the first integral equation
solver for the calculation of Taylor states in toroidal regions. While
preliminary results from this work were given in~\cite{EpGrOn2},
details of the actual solver were not provided. A separation of
variables numerical solver for the full \emph{exterior}
axisymmetric electromagnetic scattering problem from perfect
conductors is discussed in~\cite{epstein_2017}, but this work does not
address the computation of interior eigenfunctions nor solve the
boundary value problem with data on the normal components of~$\bE$, $\bH$.
We close this gap with
the present work. The integral formulation we present here applies to
both toroidally axisymmetric and non-axisymmetric domains, but thus
far our numerical solver can only treat the first situation. We will
therefore restrict the description of the numerical solver to that
case. Let us stress again that while the domain is axisymmetric, the
Taylor state itself may not be, and the solver we present here can
compute these non-axisymmetric equilibria. As such, it may be applied
to the computation of magnetohydrodynamic equilibria in
spheromaks~\cite{geddes}, reversed field
pinches~\cite{taylor_revisited}, and in tokamaks for start-up
scenarios~\cite{battaglia} and the study of
magnetohydrodynamic instabilities~\cite{hastie,gimblett}.

\begin{figure}[t]
    \centering
    \begin{subfigure}[t]{.5\linewidth}
        \centering
        \includegraphics[width=.95\linewidth]{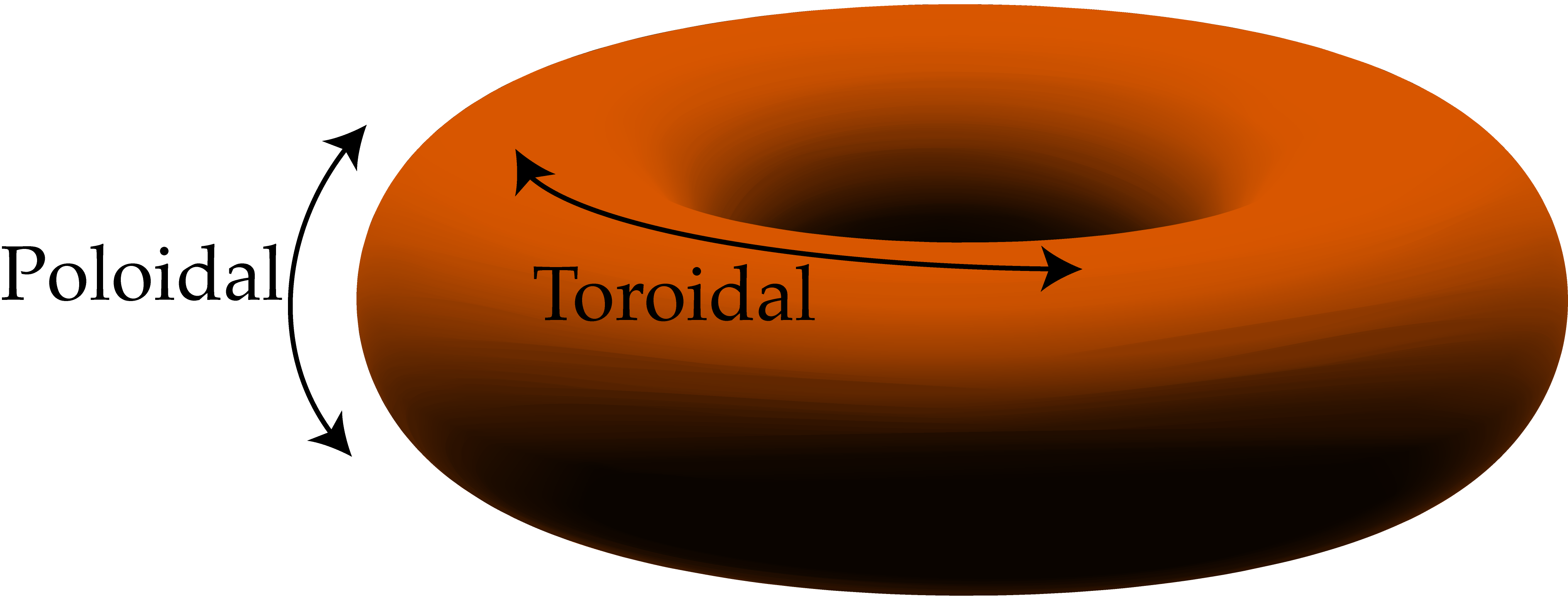}
        \caption{Poloidal and toroidal directions.}\label{fig_poloidal}
    \end{subfigure}
    \hfill
    \begin{subfigure}[t]{.4\linewidth}
        \centering
        \includegraphics[width=.95\linewidth]{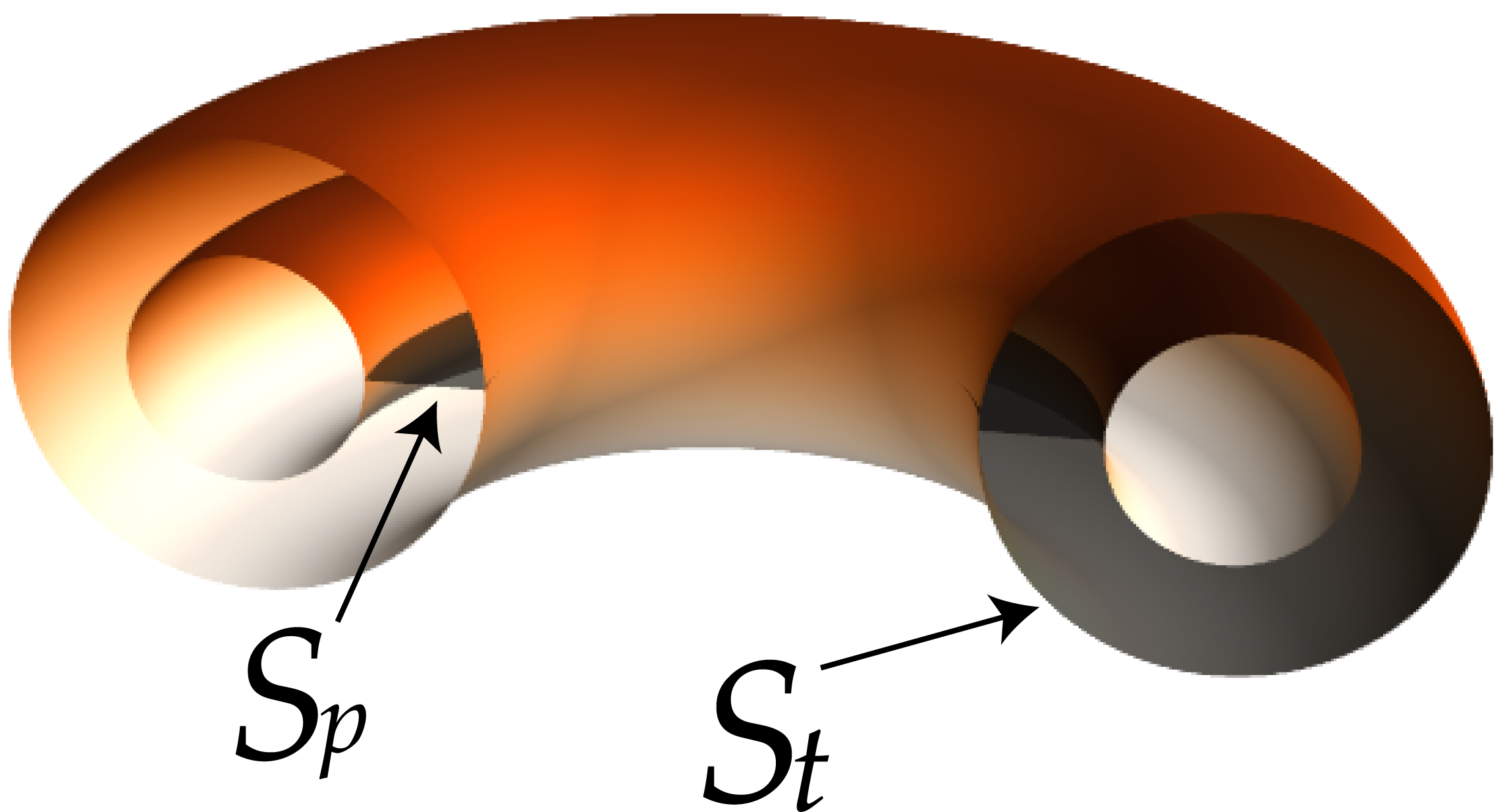}
        \caption{Genus-two toroidal shell.}\label{fig_shell}
    \end{subfigure}
    \caption{The basic geometry of the toroidal domains relevant to
      axisymmetric magnetohydrodynamics.  The boundary $\Gamma$ of the
      domain is shown in orange shading, while the surfaces $S_{t}$
      and $S_{p}$ on which the magnetic fluxes are specified are shown
      in gray shading.}
\end{figure}

The mathematically well-posed form of the problem is as follows. We
construct numerical solutions to the Beltrami boundary-value problem
given by:
\begin{equation}\label{eq-bv}
\begin{aligned}
\nabla \times \bB &= \lambda \bB & \quad & \text{in } \Omega, \\
\bB \cdot \vct{n} &= 0 & & \text{on } \Gamma,
\end{aligned}
\end{equation}
where
$\lambda$ is a real number given as input to the solver,
$\Omega$ is an axisymmetric toroidal domain, and
$\Gamma = \partial \Omega$ is the (smooth)
boundary of the region~$\Omega$. 
Depending on the genus of $\Gamma$, additional (topological)
constraints on $\bB$ must be added in order for~(\ref{eq-bv}) to be
well-posed. For Taylor states in laboratory plasmas, it is often
natural to take these constraints as conditions on the toroidal and
poloidal flux of $\bB$~\cite{taylor_review,hudson}, see
Figure~\ref{fig_poloidal}. In a genus-two toroidal flux shell (see
Figure~\ref{fig_shell}) two conditions must be imposed:
\begin{equation} \label{eq-flux}
  \int_{S_t} \bB \cdot d\ba = \fluxtor \qquad \text{and} \qquad
  \int_{S_p} \bB \cdot d\ba = \fluxpol,
\end{equation}
where $d\ba = \vct{n} \, da$, with $da$ being the surface area element
and $\vct{n}$ the oriented normal along the surfaces~$S_t$ and~$S_p$.
On the other hand, if the toroidal domain is not hollow (genus-one),
only one additional flux condition is necessary.  The need for extra
conditions~(\ref{eq-flux}) to ensure well-posedness stems from the
multiply-connectedness of the boundary $\Gamma$ -- namely, the
existence of {\em harmonic surface vector fields} on $\Gamma$ and
interior volume $\lambda$-Neumann vector fields in
$\Omega$~\cite{EpGr}. Readers interested in more details on the
well-posedness of the boundary value
problem~(\ref{eq-bv},~\ref{eq-flux}) may read
references~\cite{kress_1986,EpGrOn2}.

Our integral equation
formulation is based on the observation that if $\nabla \times \bB =
\lambda \bB$, then the pair $\{\bE,\bH\} =\{ i\bB,\bB\}$ satisfies the
time-harmonic Maxwell's equations in vacuum, with $\lambda$
playing the role of a wavenumber. Boundary conditions on the normal
component of~$\bB$ then correspond to 
boundary conditions on the normal components of~$\bE$ and $\bH$.
This fact, coupled with the symmetry of~$\bE$ and $\bH$, makes it
natural to represent~$\bB$ using generalized Debye sources, as
in~\cite{EpGr,EpGrOn2}.
Application of the boundary condition in~\eqref{eq-bv}
and flux constraints~\eqref{eq-flux} to the generalized Debye source
representation immediately yields a second-kind integral equation
which can be solved with standard techniques.

The paper is organized as follows. In Section~\ref{sec_debye}, we
establish the link between linear Beltrami fields and the generalized
Debye representation at the heart of our integral equation
formulation. In Sections~\ref{sec_toroidal} and~\ref{sec_shell}, we
derive the second-kind integral equations for the densities of the
vector and scalar potentials in the generalized Debye representation
of the Beltrami field. Section~\ref{sec_toroidal} applies to toroidal
regions, while Section~\ref{sec_shell} applies to hollow toroidal
shells. Section~\ref{sec_numerical} describes our numerical method to
compute the solution to the integral equations, and to subsequently
evaluate the
Beltrami fields. In Section~\ref{sec_examples}, we illustrate the flexibility and accuracy of our solver with three examples that
are relevant to laboratory plasma experiments, and in
Section~\ref{sec_conclusion} we summarize our work and discuss ideas
for further development.

\section{Beltrami fields and the generalized Debye source representation}\label{sec_debye}

In this section, we introduce the generalized Debye source
representation of linear Beltrami fields. The generalized Debye source
representation was originally developed for solving time-harmonic
electromagnetic scattering problems from perfect electric conductors
and dielectric materials~\cite{EpGr,EpGrOn}, but has recently been
found to be extremely well-suited for describing force-free fields
with spatially constant Beltrami parameter
$\lambda$~\cite{EpGrOn2}. The representation immediately leads to a
well-conditioned (away from {\em physical} interior resonances)
second-kind integral equation which can be numerically inverted to
high-precision. We first present the generalized Debye representation
for the time-harmonic Maxwell equations, and then make the connection
with linear Beltrami fields.

\subsection{Debye source representation for time-harmonic Maxwell's equations}\label{sec:Debye_summary}

In their simplest form, the time-harmonic Maxwell's equations in a region free of charge and
current, are given by:
\begin{equation}\label{Maxwell_timeharmonic}
  \begin{aligned}
    \nabla \times \bE &= i k \bH \qquad & 
    \nabla \times \bH &= -i k \bE \\
    \nabla \cdot \bE &= 0 & \nabla \cdot \bH &= 0,
  \end{aligned}
\end{equation}
where $\bE$ and $\bH$ are the electric and magnetic fields,
respectively, and have been suitably scaled by the electric
permittivity and magnetic permeability so that the equations only
depend on a single parameter, $k$. The real wavenumber $k$ is
proportional to the time-harmonic angular frequency $\omega$, with $k
= \omega/2\pi$. Epstein and Greengard have recently developed a new
and robust integral representation for the solution
to~\eqref{Maxwell_timeharmonic} in the exterior of objects with smooth
boundaries for 
the standard perfect electric conductor homogeneous
boundary conditions on~$\bE$ and $\bH$~\cite{EpGr}:
\begin{equation}
  \vct{n} \times \bE = 0, \qquad \vct{n} \cdot \bH = 0.
\end{equation}
This representation, called the generalized Debye source
representation, is a full potential-antipotential formulation, given
by:
\begin{equation} 
  \begin{aligned}
    \bE &= ik \bA - \nabla u - \nabla \times \bQ, \\
    \bH &= ik \bQ - \nabla v + \nabla \times \bA,
  \end{aligned}
\end{equation}
along with the two consistency conditions on the potentials
\begin{equation}\label{eq_consist}
\nabla \cdot \bA = iku, \qquad \nabla \cdot \bQ = ikv.
\end{equation}
In the exterior of a scatterer with smooth closed boundary $\Gamma$, 
the vector and scalar potentials are constructed as
\begin{equation}
\begin{aligned}
\bA(\bx) &= \int_{\Gamma} \frac{e^{ik|\bx-\bx'|}}{4\pi |\bx-\bx'|}
\, \bj(\bx') \, da(\bx'), & \qquad 
\bQ(\bx) &= \int_{\Gamma} \frac{e^{ik|\bx-\bx'|}}{4\pi |\bx-\bx'|}
\, \bmm(\bx') \, da(\bx'), \\
u(\bx) &= \int_{\Gamma} \frac{e^{ik|\bx-\bx'|}}{4\pi |\bx-\bx'|}
\, \rho(\bx') \, da(\bx'),  & 
v(\bx) &= \int_{\Gamma} \frac{e^{ik|\bx-\bx'|}}{4\pi |\bx-\bx'|}
\, \sigma(\bx') \, da(\bx'),
\end{aligned}
\end{equation}
where $da(\bx')$ denotes the surface area differential along $\Gamma$ 
at the point $\bx'$, and
the kernel is the Green's function for the Helmholtz equation in three
dimensions (Chapter 2 of~\cite{nedelec}):
\begin{equation}
g(\bx,\bx') = \frac{e^{ik\vert \bx - \bx'\vert}}{4\pi \vert \bx-\bx' \vert}.
\end{equation}
In order for the consistency conditions~\eqref{eq_consist} to be
satisfied, it must be the case that~\cite{EpGr}
\begin{equation}
\nabla_\Gamma \cdot \bj = ik\rho, \quad \text{and} \quad
\nabla_\Gamma \cdot \bmm = ik\sigma,
\label{eq:consistency_conditions}
\end{equation}
where $\nabla_\Gamma \cdot \bF$ denotes the intrinsic {\em surface
  divergence} of a tangent vector field $\bF$ along $\Gamma$.
When the scattering object is a perfect conductor,  the surface vector fields $\bj$, $\bmm$ are written as
\begin{equation}\label{eq_jm}
  \begin{aligned}
    \bj &= ik\left( \nabla_\Gamma \surflap^{-1} \rho - \vct{n}
      \times \nabla_\Gamma \surflap^{-1} \sigma \right) + \bj_H \\
    \bmm &= ik\left( \nabla_\Gamma \surflap^{-1} \sigma + \vct{n}
      \times \nabla_\Gamma \surflap^{-1} \rho \right) 
    + \bmm_H \\
    &= \vct{n} \times \bj,
  \end{aligned}
\end{equation}
and automatically satisfy \eqref{eq:consistency_conditions}. Here $\nabla_\Gamma$ is the surface gradient operator, and
by $\surflap^{-1}$ we denote the
inverse of the surface Laplacian (Laplace-Beltrami operator) along
$\Gamma$ restricted to the class of mean-zero functions. That is to
say, the generalized Debye sources $\rho$, $\sigma$ are required to
have mean zero, as well as the functions \mbox{$\alpha =
\surflap^{-1} \rho$} and $\beta =
\surflap^{-1} \sigma $. Restricted to the class of mean-zero
functions defined on $\Gamma$, the Laplace-Beltrami operator is
uniquely invertible~\cite{EpGr}.
  Furthermore, the surface vector fields $\bj_H$ and
$\bmm_H$ are {\em harmonic}, in the sense that
\begin{equation}
\nabla_\Gamma \cdot \bj_H = 0, \qquad \nabla_\Gamma \cdot \vct{n} \times \bj_H = 0,
\end{equation}
and likewise for $\bmm_H$. 
The construction of $\bj$, $\bmm$ in~\eqref{eq_jm} ensures
uniqueness of the underlying representation. 

This concludes the overview of the generalized Debye source
representation that we will use in the remainder of the paper. More
details regarding the simulation of time-harmonic Maxwell fields using
the generalized Debye source approach can be found in~\cite{EpGr,
  EpGrOn, epstein_2017}. We are now ready to turn to the problem of
interest, namely the construction of an integral representation for
linear Beltrami fields based on the generalized Debye
representation. We start by highlighting the link between the
scattering problem discussed above and linear Beltrami fields.

\subsection{Generalized Debye representation for linear Beltrami fields}

We want to compute, for a given $\lambda \geq 0$,
the magnetic field $\bB$ satisfying
\begin{equation}\label{eq-bv-repeated}
\begin{aligned}
\nabla \times \bB &= \lambda \bB & \quad & \text{in } \Omega, \\
\bB \cdot \vct{n} &= 0 & & \text{on } \Gamma,
\end{aligned}
\end{equation}
where $\Omega$ is an axisymmetric toroidal domain, and $\Gamma =
\partial \Omega$ is the boundary of $\Omega$. Note that as stated
here, \eqref{eq-bv-repeated} does not have a unique solution. For
uniqueness, additional constraints must be imposed on $\bB$, which we
will take to be constraints on the magnetic flux, as already discussed
in the introduction. We will return to this point in Sections
\ref{sec_toroidal} and \ref{sec_shell}, where we construct the
integral equations used by our solver.

Now, assume that $\bB$ satisfies \eqref{eq-bv-repeated}. Then the pair
$\{\bE$,$\bH\}$ given by
\begin{equation} \bE = i \bB, \quad
  \text{and} \quad \bH = \bB.
\end{equation}
satisfies Maxwell's equations \eqref{Maxwell_timeharmonic}, with
$\lambda$ playing the role of $k$ and the boundary conditions \mbox{$\bE
\cdot \vct{n}=0$} and \mbox{$\bH \cdot \vct{n}=0$} on
$\Gamma$. Analysis of this
interior boundary value problem was done in 1986 by Kress~\cite{kress1},
and later discussed with regard generalized Debye sources and the
Beltrami problem in 2015~\cite{EpGrOn2}, but no detailed computations were
performed.
With this in mind, motivated by the specific
relationship $\bE=i\bH$, the symmetry of the generalized Debye
representation provides a natural representation for the Maxwell pair
$\{\bE$,$\bH\}$: the vector and scalar potentials must satisfy
\cite{EpGrOn}
\begin{equation}\label{eq_iaq}
\bA = i \bQ \quad \text{and} \quad u = i v.
\end{equation}
The relationship above is required for Beltrami fields in the interior
or exterior of $\Gamma$, otherwise \mbox{$\nabla\times\bB=\lambda\bB$}
is not satisfied. We thus represent our force-free field $\bB$ as
\begin{equation}
  \bB = i\lambda \bQ - \nabla v + i \, \nabla \times \bQ,
\end{equation}
where $\bQ$ and $v$ are generalized Debye potentials for boundary
conditions associated with the magnetic field along a perfect
conductor.

We now turn our attention to the particular representation that we
will use for Beltrami fields in the interior of genus-one tori and
genus-two toroidal shells, i.e. the exact construction of the surface
vector field $\bmm$ introduced above, following the presentation given
in Section \ref{sec:Debye_summary}. The expressions for the case of
genus-one tori and of genus-two toroidal shells differ slightly on
several occasions, so we treat each case in a separate section.

\section{Beltrami fields in toroidal geometries}
\label{sec_toroidal}

In the previous section, the generalized Debye source representation
was shown to be a natural way to represent Beltrami fields because of
the inherent symmetry between the electric and magnetic fields.  
In this section, we
provide the exact construction of the surface vector field
$\bmm$ as well as derive a second-kind integral equation which is
uniquely invertible {\em except} when the Beltrami parameter~$\lambda$
is precisely an eigenvalue of the curl operator~\cite{EpGrOn2}.

As explained in the previous section, in the interior of a
genus-one torus we will represent Beltrami fields as
\begin{equation} \label{eq-repr}
\bB = i\lambda \bQ - \nabla v + i \,  \nabla \times \bQ,
\end{equation}
where $\bQ$ and $v$ are vector and scalar potentials, respectively,
given by layer potentials along $\Gamma$, the boundary of a 
toroidal region $\Omega$:
\begin{equation}\label{eq-pot_convol}
    \bQ(\bx) = \int_{\Gamma} g(\bx,\bx')
    \, \bmm(\bx') \, da(\bx'), \qquad 
    v(\bx) = \int_{\Gamma} g(\bx,\bx')
    \, \sigma(\bx') \, da(\bx'). 
\end{equation}
In order for $\bB$ to satisfy Beltrami's equation (and for the pair
$\{i\bB$, $\bB\}$ to satisfy Maxwell's equations) it is necessary that
\begin{equation}\label{eq_mconsist}
\nabla_\Gamma \cdot \bmm = i\lambda \sigma.
\end{equation}
We now provide an explicit construction of $\bmm$ in terms of $\sigma$
and a harmonic surface vector field $\bmm_H$ such that the
representation is unique and leads to an invertible integral
equation. 
To this end, let us assume that all differential surface
operators are oriented with respect to a unit normal vector $\vct{n}$
along $\Gamma$
which is always assumed to point into the region $\Omega^c$, the
complement of $\Omega$
(i.e. the unit {\em outward} normal).  
As described in Section \ref{sec:Debye_summary}, the surface vector field 
$\bmm$ given by
\begin{equation} \label{eq-mbuild}
\bmm = i\lambda \left( \nabla_\Gamma \surflap^{-1} \sigma -
i\, \vct{n} \times \nabla_\Gamma \surflap^{-1} \sigma \right) + \alpha \,
\bmm_H.
\end{equation}
is constructed to automatically enforce the consistency
condition \eqref{eq_mconsist} and yield uniqueness in the resulting
integral equation, derived below, in~\eqref{eq-inteq1}.
The surface vector field
$\bmm_H$ appearing in expression~\eqref{eq-mbuild} is a tangential
 harmonic vector field satisfying the conditions:
\begin{equation}\label{eq-surfharm}
\begin{gathered} 
\nabla_{\Gamma} \cdot \bmm_H = 0, \qquad \quad
\nabla_{\Gamma} \cdot \vct{n} \times \bmm_H = 0, \\
\vct{n} \times \bmm_H = i \, \bmm_H.
\end{gathered}
\end{equation}
The constant $\alpha$ appearing in front of $\bmm_H$ is a complex
number to be determined in the solution of the Beltrami boundary-value
problem. The relationship between~$\vct{n} \times \bmm_H$ and $\bmm_H$ is
in fact necessary in order to ensure uniqueness of the integral
representation~\cite{EpGrOn2}, and is analogous to the relationship
between~$\bA$ and~$\bQ$ in~\eqref{eq_iaq}.

Using this representation for $\bB$, it is straightforward to derive a
second-kind integral equation (augmented with flux conditions) for the
unknowns $\sigma$ and $\alpha$ by merely enforcing the local boundary
condition $\bB \cdot \vct{n} = 0$ and global integral
constraint. Introducing the {\em single-layer operator} $\cS$ as
\begin{equation}
\cS[ f](\bx) = \int_{\Gamma} g(\bx,\bx') \, f(\bx') \, da(\bx'),
\end{equation}
the integral equation and flux condition can be written as
\begin{equation}\label{eq-inteq1}
\begin{aligned}
\frac{\sigma}{2}+ \cK [\sigma, \alpha] &= 0, \\
\int_{S_t} \bB(\sigma,\alpha) \cdot d\ba &= \fluxtor,
\end{aligned}
\end{equation}
where $d\ba = \vct{n} \, da$ and with the
notation~$\bB(\sigma,\alpha)$ we highlight the fact that the magnetic
field~$\bB$ is a function of~$\sigma$ and~$\alpha$.
In~\eqref{eq-inteq1} the operator $\cK$ has been used to abbreviate the
compact operator
\begin{equation}\label{eq_compact}
\mathcal K [\sigma,\alpha] = i\lambda \, \vct{n} \cdot \mathcal
S [\bmm] - \mathcal \cS'[\sigma]+ i \, \vct{n} \cdot \nabla \times \mathcal S [\bmm],
\end{equation}
where $\bmm = \bmm(\sigma,\alpha)$ as in~\eqref{eq-mbuild}, and $\cS'$ is the normal derivative of the single-layer
potential, interpreted
in the appropriate principal-value sense~(Chapter 2 of \cite{colton_kress}).
  The flux condition in
system~(\ref{eq-inteq1}) is still not optimal from the point of view
of a boundary integral equation because the integration is being
carried out over an additional surface $S_t$ (depicted in
Figure~\ref{fig_shell}) spanning the cross section
of the toroidal region $\Omega$. Letting $\partial S_t$ denote the curve
bounding this surface, $\partial S_t = \Gamma \cap S_t$, we can reduce
this flux integral to a circulation integral using Stokes' theorem
and the fact that $\bB$ is a Beltrami field:
\begin{equation}
\begin{aligned}
\int_{S_t} \bB \cdot d\ba &= \int_{S_t} \frac{1}{\lambda} (\nabla
\times \bB)
\cdot d\ba, \\ 
&= \frac{1}{\lambda} \int_{\partial S_t} \bB \cdot d\bl,
\end{aligned}
\end{equation}
where $d\bl$ denotes the unit arclength differential. 
For small values of $\lambda$,
special care must be taken in numerically evaluating this integral,
as discussed in Section~5 of \cite{EpGrOn}.

We now turn our attention to deriving an analogous representation for
Beltrami fields in genus-two toroidal shells. The only difference with the present situation is
the inclusion of an extra surface harmonic vector field because of the
genus of the boundary. For this
reason, two integral (flux) constraints must be enforced to
construct a well-posed boundary-value problem, as opposed to the single
condition enforced in this section.

\section{Beltrami fields in toroidal shells}
\label{sec_shell}

In this section, we extend the representation of Beltrami fields
introduced in the previous section to toroidal geometries with genus
two.  Specifically, we turn our attention to geometries which are
topologically equivalent to that in Figure~\ref{fig_shell}, which from
now on we will refer to as a {\em toroidal shell}. The viability of
this extension to geometries whose boundary consists of more than
one component is discussed in~\cite{EpGrOn2}.  The following
discussion makes one minor change in notation from the previous
section: the interior of the genus-two toroidal shell is still given
by $\Omega$, however its boundary $\Gamma$ is now explicitly given by
the union of two disjoint boundaries, $\Gamma = \Gamma^\myout \cup
\Gamma^\myin$. The unbounded component of $\Omega^c$ is assumed to be
bounded by $\Gamma^\myout$ and the bounded component is assumed to be
bounded by $\Gamma^\myin$.  We assume that the unit normal $\vct{n}$
along $\Gamma^\myout$ points into the unbounded component of
$\Omega^c$, but that the unit normal $\vct{n}$ along $\Gamma^\myin$
points into $\Omega$. This convention may be non-standard, but it
treats each surface as having the same counter-clockwise
parameterization and simplifies implementation details.

\begin{figure}[t]
    \centering
    \includegraphics[width=.5\linewidth]{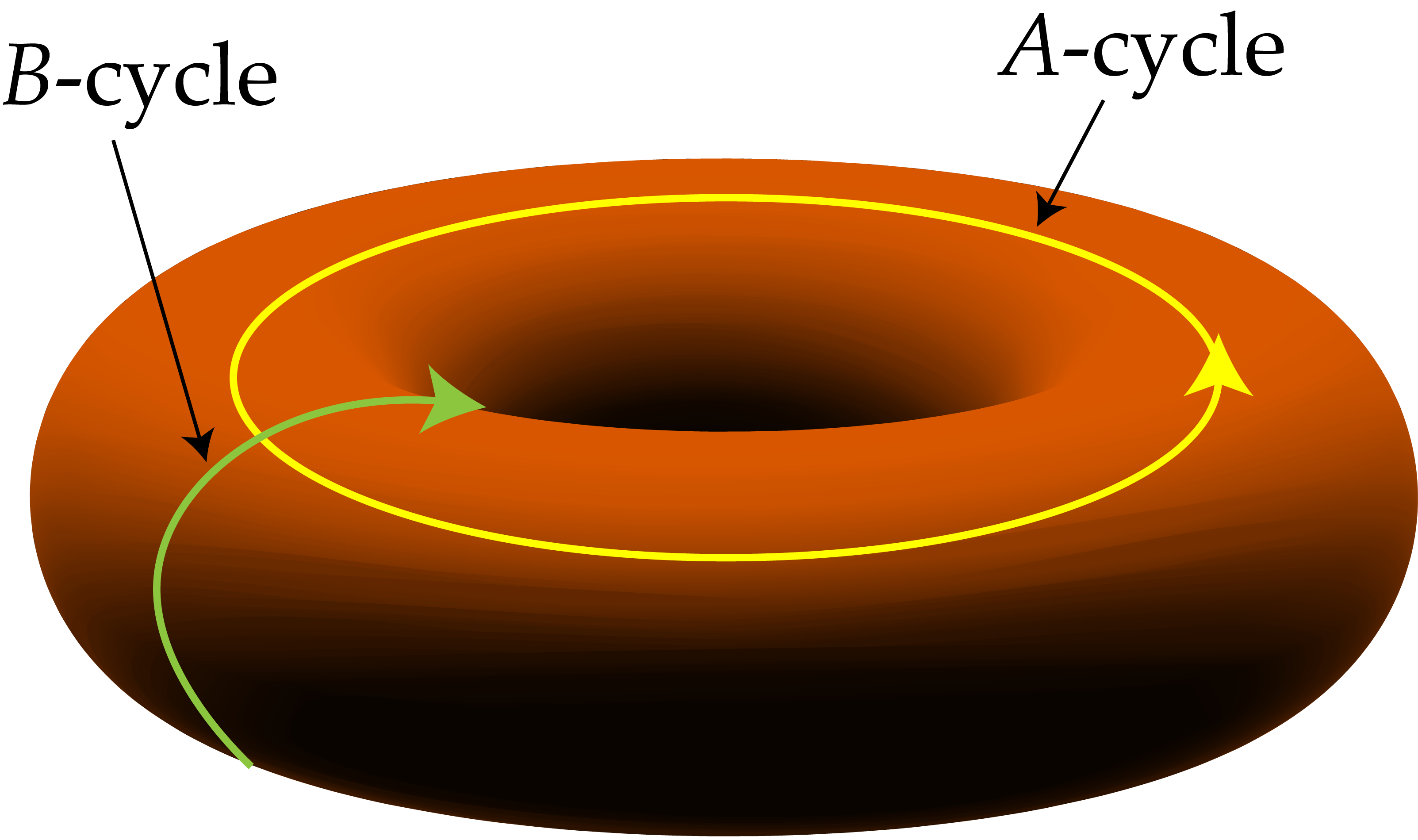}
    \caption{$A$- and $B$-cycles along a toroidal boundary.}
    \label{fig_cycles}
\end{figure}

The representation for Beltrami fields in $\Omega$ given by~\eqref{eq-repr}
does not change, however the
construction of the surface vector field $\bmm$ is altered slightly to
account for the disjoint boundaries. In the genus-two case, 
the vector field $\bmm$ consists of two components, $\bmm^\myout$
defined along $\Gamma^\myout$ and $\bmm^\myin$ defined along
$\Gamma^\myin$:
\begin{equation}
\bmm =
\begin{cases}
i\lambda\left( \nabla_\Gamma \surflap^{-1} \sigma - i\, \vct{n}
\times \nabla_\Gamma \surflap^{-1} \sigma \right) + \alpha \, \bmm^\myout_H & 
\qquad \text{on } \Gamma^\myout, \\
i\lambda\left( \nabla_\Gamma \surflap^{-1} \sigma + i\, \vct{n}
\times \nabla_\Gamma \surflap^{-1} \sigma \right) + \beta \, \bmm^\myin_H & 
\qquad \text{on } \Gamma^\myin,
\end{cases}
\end{equation}
where the harmonic surface vector fields $\bmm^\myout_{H}$ and
$\bmm^\myin_{H}$ (which are only defined on $\Gamma^\myout$ and
$\Gamma^\myin$, respectively) must satisfy:
\begin{equation}
\begin{aligned}
\vct{n} \times \bmm_H^\myout &= i  \bmm^\myout_{H}, \\ 
\vct{n} \times \bmm_H^\myin &= -i \, \bmm^\myin_{H}.
\end{aligned}
\end{equation}
The difference in these relationships is due to the fact that the
Beltrami field in $\Omega$ appears as an {\em exterior} field from the
point of view of $\Gamma^\myin$.
As in the previous section, enforcing the
boundary conditions
\begin{equation}
\begin{gathered}
\bB \cdot \vct{n} = 0 \qquad \text{on } \Gamma, \\
\int_{S_t} \bB \cdot d\ba = \fluxtor, \qquad 
\int_{S_p} \bB \cdot d\ba = \fluxpol
\end{gathered}
\end{equation}
yields an augmented second-kind integral equation for $\sigma$,
$\alpha$, and $\beta$ which is similar to the previous situation:
\begin{equation}\label{eq-inteq2}
\begin{aligned}
\frac{1}{2}\sigma^\myout 
  + \mathcal K [\sigma, \alpha, \beta] &= 0, &\quad 
&\text{on } \Gamma^\myout,\\
-\frac{1}{2}\sigma^\myin 
+ \mathcal K [\sigma, \alpha, \beta] &= 0, & &\text{on } \Gamma^\myin,\\
\int_{S_t} \bB(\sigma,\alpha,\beta) \cdot d\ba &= \fluxtor, & & \\
\int_{S_p} \bB(\sigma,\alpha,\beta) \cdot d\ba &= \fluxpol, & & 
\end{aligned}
\end{equation}
where we have once more emphasized that the magnetic field $\bB$ is
expressed in terms of $\sigma$, $\alpha$ and $\beta$ and 
where we have explicitly shown the dependence
on densities $\sigma^\myout$,~$\sigma^\myin$ which are defined along
$\Gamma^\myout$ and $\Gamma^\myin$, respectively, to highlight the
change of sign on the diagonal.  The compact operator $\mathcal K$ is
as before,
\begin{equation}
\mathcal K [\sigma,\alpha,\beta] 
= i\lambda \, \vct{n} \cdot \cS [\bmm]  - \cS'[\sigma] 
+ i \, \vct{n} \cdot \nabla \times \mathcal S [\bmm].
\end{equation}
except that we have to keep in mind that the operator has contributions from densities on both $\Gamma^\myout$ and $\Gamma^\myin$. 
The flux surface integrals in~\eqref{eq-inteq2} can be written 
as circulation
integrals on $\Gamma^\myout$ and $\Gamma^\myin$ by invoking Stokes' Theorem:
\begin{equation}
\int_{S_t} \bB \cdot d\ba = 
  \frac{1}{\lambda} \int_{\partial S_1} \bB \cdot d\bl, \qquad 
\int_{S_p} \bB \cdot d\ba =
  \frac{1}{\lambda} \int_{\partial S_2} \bB \cdot d\bl,
\end{equation}
where the line integrals are on the closed curves $\partial
S_{1}=\Gamma\cap S_{t}$ and $\partial S_{2}=\Gamma \cap
S_{t}$. Both of these flux integrals can be computed as
the difference of circulation integrals, the former around what are
known as $B$-cycles and the latter around
$A$-cycles, see Figure~\ref{fig_cycles}. We describe the
actual numerical calculation of these integrals in
Section~\ref{sec_numerical}.

\section{An axisymmetric integral equation solver}
\label{sec_numerical}

In this section we describe the numerical implementation of the integral equation solver we have developed to calculate force-free Beltrami fields in
toroidally axisymmetric geometries. A similar version of this solver,
optimized for full electromagnetic
multi-mode calculations, appears in~\cite{epstein_2017}.
Related to this work are the 
very efficient scalar axisymmetric solvers for the
Laplace and Helmholtz equations presented
in~\cite{young, helsing_2016, helsing_2015, helsing_2014}.
In the following discussion, 
we assume that a point $\bx \in \mathbb R^3$ has cylindrical
coordinates given by $\bx = (x,y,z) = (r\cos\phi,r\sin\phi,z)$, where $r$ represents the distance from the axis of revolution, and $\phi$ is the toroidal angle.
The orthonormal unit vectors in Cartesian coordinates will be denoted
$\vct{i}$, $\vct{j}$, $\vct{k}$, and the orthonormal unit vectors
in cylindrical
coordinates will be denoted in bold-face as $\vct{r}$, $\bphi$, $\vct{z}$. The Cartesian and cylindrical bases are related through
\begin{equation}
  \begin{aligned}
    \vct{r}(\phi) &= \cos\phi \, \vct{i} + \sin\phi \, \vct{j},
    &\qquad \vct{i} &= \cos\phi \, \vct{r}(\phi) - \sin\phi \,
    \bphi(\phi), \\
    \bphi(\phi) &= -\sin\phi \, \vct{i} + \cos\phi \, \vct{j},
    &\qquad \vct{j} &= \sin\phi \, \vct{r}(\phi) + \cos\phi \,
    \bphi(\phi), \\
    \vct{z} &= \vct{k}, & \vct{k} &= \vct{z},
    \end{aligned}
\end{equation}
These relationships will be useful in
the following sections. The following addition formulae will also be
useful:
\begin{equation}\label{eq_vecaddition}
  \begin{aligned}
    \vct{r}(\phi-\phi') &= \cos\phi' \, \vct{r}(\phi) -
    \sin\phi' \, \bphi(\phi) \\
    \bphi(\phi-\phi') &= \sin\phi' \, \vct{r}(\phi) +
    \cos\phi' \, \bphi(\phi) \\
  \end{aligned}
\end{equation}
The dependence of $\vct{r}$ and $\bphi$ on $\phi$ 
will be dropped when unnecessary. 

\begin{figure}[t]
    \centering
    \begin{subfigure}[t]{.5\linewidth}
        \centering
        \includegraphics[width=.95\linewidth]{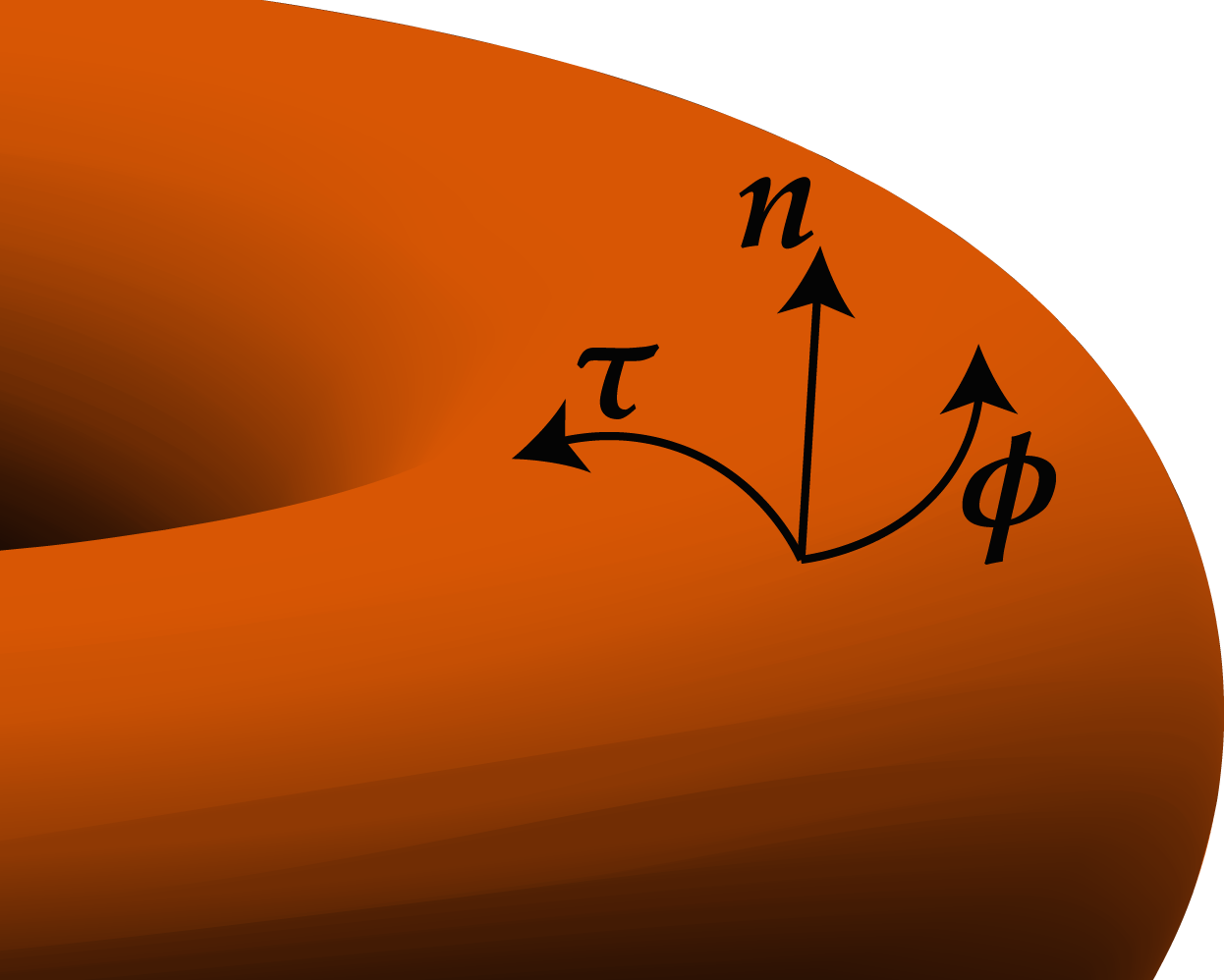}
        \caption{Local surface coordinate system along
          $\Gamma^\myout$.}
        \label{fig_coordinates}
    \end{subfigure}
    \hfill
    \begin{subfigure}[t]{.45\linewidth}
        \centering
        \includegraphics[width=.8\linewidth]{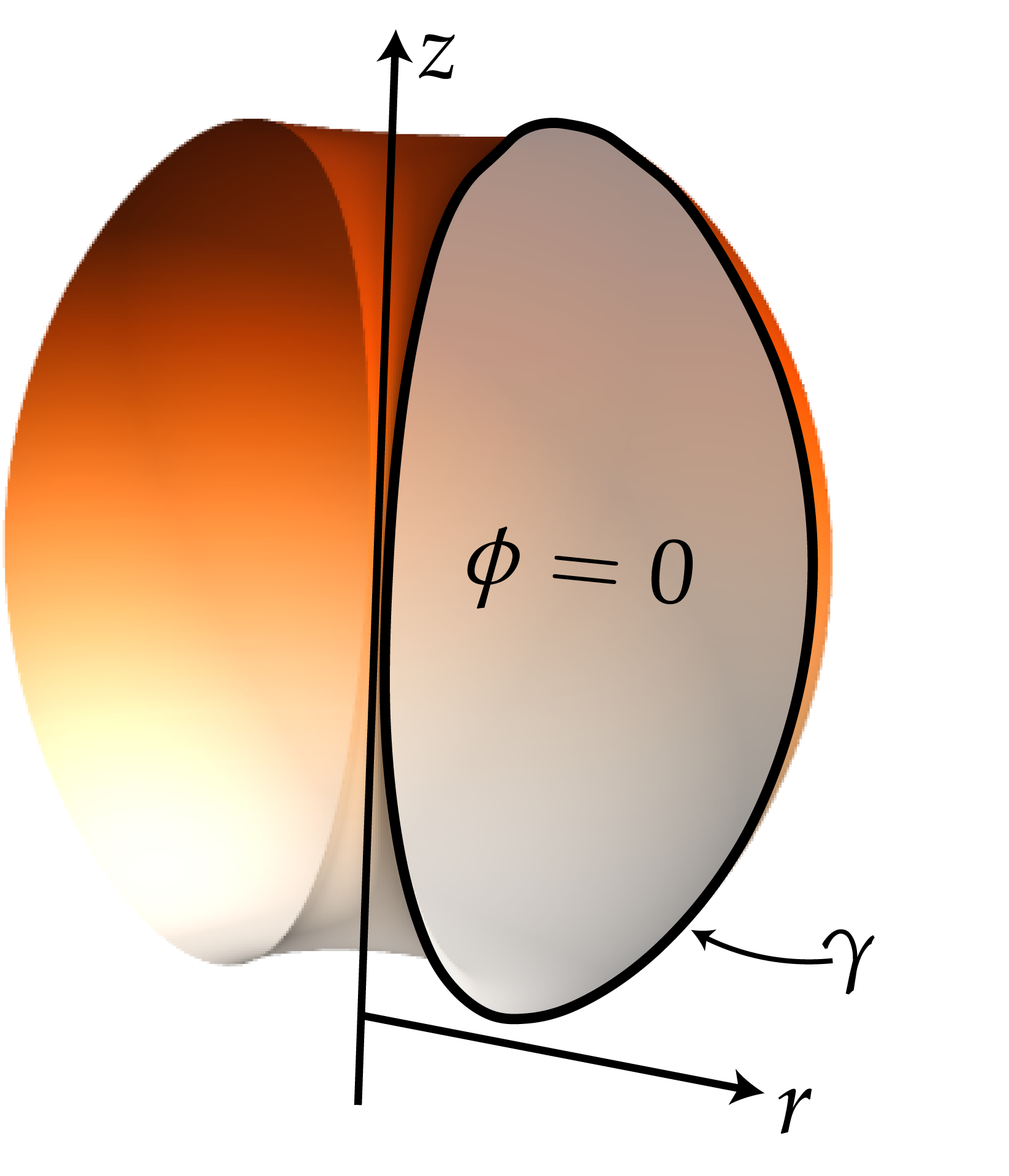}
        \caption{Generating curves for toroidal shells.}\label{fig_gen}
    \end{subfigure}
    \caption{Geometry parameterization. The shell~$\Gamma$ is shaded
      in orange, and the generating curve is drawn in black in
      Figure~\ref{fig_gen}.}
\end{figure}

We now revisit some earlier notation. Let $\Omega$, as before, be a
toroidally axisymmetric domain in~$\mathbb R^3$ with
boundary $\Gamma$. That is
to say, rotation about the $z$-axis does not alter $\Omega$.
The
boundary of the $\phi=0$ cross-section of $\Omega$ in the $xz$-plane
(and the $rz$-plane in cylindrical coordinates) will be denoted by
$\gamma$, and will be referred to as the {\em generating curve} for
$\Gamma$, see Figure~\ref{fig_gen}. Note that in the $\phi=0$ plane,
\begin{equation}
\vct{r} = \vct{i}, \quad \bphi = \vct{j}, \quad  \vct{z} = \vct{k}.
\end{equation}
 In a minor abuse of notation, we
will assume that the curve $\gamma: [0,L] \to \mathbb R^2$ is parameterized
 by functions $r(s)$, $z(s)$ of the arc length $s$:
\begin{equation}
\gamma(s) = (r(s), z(s)),
\end{equation}
and therefore $\vert d\gamma/ds \vert = 1$.
Any point $\bx \in \Gamma$ is given as
\begin{equation}
\bx(s,\phi) = r(s)  \cos\phi \, \vct{i} +
r(s)  \sin\phi \, \vct{j} + z(s) \vct{k}.
\end{equation}
Furthermore, we will assume that $\gamma$ is
positively oriented such that the unit normal along $\Gamma$, given by
\begin{equation}
\vct{n} = \frac{dz}{ds} \vct{r} -\frac{dr}{ds} \vct{z},
\end{equation}
points into $\Omega^c$. Note that $\vct{n} = \vct{n}(s,\phi)$.
Along the boundary $\Gamma$, it will be useful to work in a system of 
local coordinates $\btau$, $\bphi$, $\vct{n}$, shown in
Figure~\ref{fig_coordinates}.
Relative to the coordinate system given by $\vct{r}$, $\bphi$, $\vct{z}$, the
unit vector $\btau$ along $\Gamma$ is given by
\begin{equation}
\btau = \frac{dr}{ds} \vct{r} + \frac{dz}{ds} \vct{z}.
\end{equation}
This local coordinate system satisfies the relation $\vct{n} = \bphi
\times \btau$.

We start the description of the numerical implementation by explaining
how the surface integral equation naturally separates into Fourier
modes, as was already shown in a different context in~\cite{young},
and how the vector potential can be decomposed into Fourier modes.

\subsection{Separation into Fourier modes}
\label{sec-four}

For any integral equation on an axisymmetric surface $\Gamma$ whose
kernel $g(\bx,\bx')$ only depends on the quantity $\phi-\phi'$,
i.e. $g(\bx,\bx') = g(r,z,r',z',\phi-\phi')$ one may Fourier transform the
surface integral equation on $\Gamma$ into a series of decoupled integral
equations along its generating curve $\gamma$.
In what follows, for clarity we will use the notation:
\begin{equation}
\begin{aligned}
g(s,s',\phi-\phi') &= g(r(s),z(s),r(s'),z(s'),\phi-\phi').
\end{aligned}
\end{equation}
Since any point on $\Gamma$ is given as $\bx =
(r(s),\phi,z(s))$, the solution $\sigma$ to 
\begin{equation}
\sigma(\bx) + \int_{\Gamma} g(\bx,\bx') \, \sigma(\bx') \, da(\bx')
= f(\bx)
\end{equation}
can be written as
\begin{equation}
\sigma(s,\phi) = \sum_\ell \sigma_\ell(s) \, e^{i\ell\phi}.
\end{equation}
The functions $\sigma_\ell$  are the solution to a series of
uncoupled integral equations on $\gamma$:
\begin{equation}\label{eq_intfour}
\sigma_\ell(s) + 2\pi \int_\gamma g_\ell(s,s') \, \sigma_\ell(s') \, r(s')
\, ds'
= f_\ell(s),
\end{equation}
with the kernels $g_\ell$ and right hand sides $f_\ell$ given by:
\begin{equation}
\begin{aligned}
g_\ell(s,s') &= 
\frac{1}{2\pi} \int_{-\pi}^\pi g(s,s',\phi)
\, e^{-i\ell\phi} \, d\phi, \\
f_\ell(s) &= \frac{1}{2\pi} \int_{-\pi}^\pi f(s,\phi) \, e^{-i\ell\phi} \, d\phi.
\end{aligned}
\end{equation}
It now remains to see how one writes the Fourier mode of the vector
potential $\bQ$ given by \eqref{eq-pot_convol}. This requires
computing the Fourier modes of a layer potential applied to a
tangential vector-field, $\bu$. To
this end, consider the following integral along $\Gamma$ for the
tangent vector field $\bu$:
\begin{equation}\label{eq_vecinteq}
\int_\Gamma g(\bx,\bx') \, \bu(\bx') \, da(\bx'),
\end{equation}
where $\bu$ is given by:
\begin{equation}\label{eq_tanvecs}
\bu = u^\tau  \btau + u^\phi  \bphi
\end{equation}
Since $\Gamma$ is a surface of revolution, the tangential vector
field in~\eqref{eq_tanvecs} can be written in terms of its Fourier
series, component-wise:
\begin{equation}
    \bu(s,\phi) = \sum_\ell \left( u^\tau_\ell(s) \, \btau
    + u^\phi_\ell(s) \, \bphi \right) e^{i\ell\phi} ,
\end{equation}
 Plugging the Fourier expansion for~$\bu$ into the integral
 in~\eqref{eq_vecinteq} we have:
\begin{equation}
 \int_\Gamma g(\bx,\bx') \, \bu(\bx') \, da(\bx') = \int_\Gamma
 g(s,s',\phi-\phi') \sum_\ell \left( u^\tau_\ell(s') \,
 \btau(s',\phi') + u^\phi_\ell(s') \, \bphi(\phi') \right)
 e^{i\ell\phi'} \, da(s',\phi'),
\end{equation}
where $da(s',\phi') = r(s') \, ds' \, d\phi'$.
 We now
calculate the projection of this integral onto the $\ell^\text{th}$
Fourier mode with respect to the unit vectors $\vct{r}(\phi)$,
$\bphi(\phi)$ at the {\em target} point $\bx$.  
The unit vectors $\vct{r}(\phi')$ and
$\bphi(\phi')$ {\em cannot} be pulled outside of the integral. However,
as a direct consequence of the addition formulas in
equation~\eqref{eq_vecaddition}, we have that the $d\phi'$ integral
can be computed as:
\begin{equation}\label{eq_intadd1}
  \begin{aligned}
    \int_{-\pi}^\pi g(\cdot,\phi-\phi') \, \vct{r}(\phi') \,
    e^{i\ell\phi'} \, d\phi' &= e^{i\ell\phi}
    \int_{-\pi}^\pi g(\cdot,\theta) \, \vct{r}(\phi-\theta) \,
    e^{-i\ell\theta} \, d\theta \\
    &= e^{i\ell\phi} \, \vct{r}(\phi) \int_{-\pi}^\pi g(\cdot,\theta) \,
    \cos\theta \,  
    e^{-i\ell\theta} \, d\theta  \\
    &\qquad \qquad -
    e^{i\ell\phi} \, \bphi(\phi) \int_{-\pi}^\pi g(\cdot,\theta) \,
    \sin\theta  \, 
    e^{-i\ell\theta} \, d\theta  \\
    &= 2\pi \, e^{i\ell\phi} \left(  g_\ell^{\cos}(\cdot) \, \vct{r}(\phi) -
        g_\ell^{\sin}(\cdot) \, \bphi(\phi) \right) 
  \end{aligned}
\end{equation}
where
\begin{equation}
  \begin{aligned}
  g_\ell^{\cos}(\cdot) &= \frac{1}{2\pi}
  \int_{-\pi}^\pi g(\cdot,\theta) \, \cos\theta
  \, e^{-i\ell\theta} \, d\theta, \\
  g_\ell^{\sin}(\cdot) &=  \frac{1}{2\pi}
  \int_{-\pi}^\pi g(\cdot,\theta) \, \sin\theta
  \, e^{-i\ell\theta} \, d\theta,
\end{aligned}
\end{equation}
and likewise
\begin{equation}\label{eq_intadd2}
\int_{-\pi}^\pi g(\cdot,\phi-\phi') \, \bphi(\phi') \,
e^{i\ell\phi'} \, d\phi' =
2\pi \, e^{i\ell\phi} \left(  g_\ell^{\sin}(\cdot) \, \vct{r}(\phi) +
        g_\ell^{\cos}(\cdot) \, \bphi(\phi) \right).
\end{equation}
We conclude that in our formulation, we do not only need to evaluate the kernel $g_{\ell}$, but we also need to compute $g_\ell^{\cos}$ and $g_\ell^{\sin}$. We address these questions in Section \ref{sec:kernel}.


Now, when solving for force-free 
Beltrami fields using~\eqref{eq-inteq1} or~\eqref{eq-inteq2}, the Fourier decomposition we use here can simplify dramatically. Indeed, we first note that for $\ell \neq 0$, any smooth
vector field $\bF$ in the interior of a torus or toroidal shell of
the form $\bF(r,z) \, e^{i\ell\phi}$ does {\em not} contribute to
toroidal flux {\em nor} poloidal flux. This can be shown via
a straightforward circulation integral, 
or by the exactness of functions
of this form (Section~7.3 of~\cite{EpGrOn2}). In other words, only
toroidally axisymmetric solutions~$\bB$ make non-zero contributions to the
flux integrals in~\eqref{eq-inteq1} and~\eqref{eq-inteq2}. 
Second, we point out that the right hand side of integral
equations~\eqref{eq-inteq1} and~\eqref{eq-inteq2} is zero except for
the flux conditions. Unless the Beltrami parameter~$\lambda$ is a
resonant value for some non-zero Fourier mode, these two facts mean
that all solutions to these integral equations have no 
angular $\phi$ dependence.
This means that we only need to compute the solution $\sigma_0$, which
satisfies the integral equation:
\begin{equation}\label{eq_inteq0}
\begin{gathered}
\frac{1}{2} \sigma_0 + \mathcal
K_0[\sigma_0,\cdot] = 0, \\
+\text{ suitable flux constraints},
\end{gathered}
\end{equation}
with~$\cK_0$ given by
\begin{equation}
  \cK_0[\sigma_0,\alpha] = i\lambda \vct{n} \cdot \cS_0 [\bmm_0] -
  \cS'_0 [\sigma_0] + i \vct{n} \cdot \nabla \times \cS_0[\bmm_0].
\end{equation}
The term $\cS_0[\bmm_0]$ can be computed using the addition
formulas in~\eqref{eq_intadd1} and~\eqref{eq_intadd2}.

In the case where $\lambda$ is a resonance corresponding to an
eigenfunction with non-trivial azimuthal dependence (i.e. an
azimuthal dependence of $e^{i\ell\phi}$ with $\ell>0$), equations
corresponding to higher azimuthal modes may need to be
solved~\cite{EpGrOn2}.  The formulae for general $\ell$ given in \eqref{eq_intfour}--\eqref{eq_intadd2} are then required, as done
in~\cite{epstein_2017} in a different context.

\subsection{Evaluation of kernels}\label{sec:kernel}

From the previous section we conclude that in order to implement a
numerical solver for our two-dimensional integral equations, it is
necessary that certain {\em modal Green's functions} be
computed~\cite{conway_cohl, helsing_2014, young}, namely the kernels
$g_\ell$, $g^{\cos}_\ell$, and $g^{\sin}_\ell$.  In particular, for
purely axisymmetric solutions, $\ell = 0$, we need to compute $g_0$,
$g^{\cos}_0$, and $g^{\sin}_0$.  The kernel appearing in Beltrami
integral equations is:
\begin{equation}\label{eq_kernel}
\begin{aligned}
g(\bx,\bx') &=
\frac{e^{i\lambda\vert\bx-\bx'\vert}}{4\pi\vert\bx-\bx'\vert} \\
&= \frac{e^{i\lambda
      R(\theta)}}{4\pi R(\theta)},
\end{aligned}
\end{equation}
where for fixed $r$, $z$, $r'$, $z'$, and $\theta = \phi-\phi'$,
the cylindrical distance $R(\theta)$ is given by
\begin{equation}
R(\theta) = \sqrt{r^2 + {r'}^2 - 2rr' \cos\theta + (z-z')^2}.
\end{equation}
We therefore now describe a method for computing:
\begin{equation}\label{eq-kern0}
\begin{aligned}
  g_0(r,r',z,z') &= \frac{1}{2\pi} \int_{-\pi}^\pi \frac{e^{i\lambda
      R(\theta)}}{4\pi R(\theta)} \, d\theta, \\
  g^{\cos}_0(r,r',z,z') &= \frac{1}{2\pi} \int_{-\pi}^\pi \frac{e^{i\lambda
      R(\theta)}}{4\pi R(\theta)} \, \cos\theta \, d\theta, \\
  g^{\sin}_0(r,r',z,z') &= \frac{1}{2\pi} \int_{-\pi}^\pi \frac{e^{i\lambda
      R(\theta)}}{4\pi R(\theta)} \, \sin\theta \, d\theta.
\end{aligned}
\end{equation}
Obviously, the integrands in~(\ref{eq-kern0}) are singular when $R=0$,
or rather when $r=r'$, $z=z'$, and $\theta=0$. Direct integration
via adaptive Gaussian quadrature, although computationally expensive,
 provides nearly full machine-precision
evaluation of $g_0$, provided $R$ is not too small and $\lambda$ is not
too large (the Beltrami parameter~$\lambda$ affects 
the number of oscillations of the integrand $[0,2\pi]$).
To accelerate the integration slightly, $g_0$ can be rewritten as
\begin{equation}
\begin{aligned}
  g_0(r,r',z,z') &= \frac{1}{2\pi} \int_{-\pi}^\pi \frac{e^{i\lambda
      R(\theta)}-1}{4\pi R(\theta)} \, d\theta + \frac{1}{2\pi}
  \int_{-\pi}^\pi \frac{1}{4\pi R(\theta)} \, d\theta \\ 
&= \frac{1}{2\pi}
  \int_{-\pi}^\pi \frac{e^{i\lambda R(\theta)}-1}{4\pi R(\theta)} \, d\theta +
  \frac{1}{4\pi^2\sqrt{rr'}} \, Q_{-\frac{1}{2}}(\Chi) \\
&= \frac{1}{2\pi}
  \int_{-\pi}^\pi \frac{e^{i\lambda R(\theta)}-1}{4\pi R(\theta)} \, d\theta +
  \frac{1}{4\pi^2\sqrt{rr'}} \sqrt{\frac{2}{1+\Chi}} 
\, K \left( \sqrt{\frac{2}{1+\Chi}}\right),
\end{aligned}
\end{equation}
where $Q_{-\frac{1}{2}}$ is the Legendre function of the second-kind of
degree negative one-half, $K$ is the complete elliptic integral of the first
kind,
\begin{equation}
K(t) = \int_0^{\frac{\pi}{2}} \frac{d\theta}
{\sqrt{1-t^2 \sin^2 \theta}} ,
\end{equation}
and the argument $\Chi$ is given by
\begin{equation}
\Chi = \frac{r^2 +{r'}^2 + (z-z')^2}{2rr'}
\end{equation}
See~\cite{cohl_1999} for a derivation of Fourier modes
of the Laplace kernel in cylindrical coordinates. The elliptic integral $K$ can be evaluated to arbitrary precision via
a combination of adaptive integration and asymptotic series near
$\theta=1$, or via Carlson's Algorithm~\cite{gil_2007,carlson_1979}. For very small $\theta$, the quantity $(e^{i\lambda R(\theta)} -
1)/R(\theta)$ can be computed to high accuracy via a Taylor expansion.
More sophisticated methods for the numerical evaluation of the function $g_\ell$, for all
$\ell$, are discussed in~\cite{helsing,young,epstein_2017}.
We compute the necessary partial derivatives of $g_0$ via analytical
differentiation and direct adaptive Gaussian quadrature, using the
fact that~\cite{nist}
\begin{equation}
 K'(t) = \frac{E(t)
  - (1 - t^2) K(t)}{t (1-t^2)},
\end{equation}
where $E(t)$ denotes the complete elliptic integral:
\begin{equation}
\begin{aligned}
E(t) &= \int_0^{\frac{\pi}{2}} \sqrt{1-t^2 \sin^2 \theta} \, d\theta, \\
E'(t) &= \frac{E(t) - K(t)}{t}.
\end{aligned}
\end{equation}

With regard to the kernels $g_0^{\cos}$ and $g_0^{\sin}$, we
first observe that:
\begin{equation}
g_\ell^{\cos} = \frac{1}{2} \left( g_{\ell+1} + g_{\ell-1}\right),
\qquad
g_\ell^{\cos} = \frac{1}{2i} \left( g_{\ell+1} - g_{\ell-1}\right).
\end{equation}
Therefore, we merely need to compute $g_{1}$ and $g_{-1}$ in order to
compute $g_0^{\cos}$ and $g_0^{\sin}$. Furthermore, since the
kernel~$g$ in~\eqref{eq_kernel} is an even function, $g_1 = g_{-1}$.
The kernel $g_1$ can be computed similarly to $g_0$ as above:
\begin{equation}
\begin{aligned}
  g_1(r,r',z,z') &= \frac{1}{2\pi} \int_{-\pi}^\pi \frac{e^{i\lambda
      R(\theta)}-1}{4\pi R(\theta)} e^{i\theta} \, d\theta + \frac{1}{2\pi}
  \int_{-\pi}^\pi \frac{e^{i\theta} }{4\pi R(\theta)} \, d\theta \\ 
&= \frac{1}{2\pi}
  \int_{-\pi}^\pi \frac{e^{i\lambda R(\theta)}-1}{4\pi R(\theta)}
  e^{i\theta} \, d\theta +
  \frac{1}{4\pi^2\sqrt{rr'}} \, Q_{\frac{1}{2}}(\Chi),
\end{aligned}
\end{equation}
where $Q_{\frac{1}{2}}$ is given by:
\begin{equation}
Q_{\frac{1}{2}}(\Chi) = \frac{1}{4\pi^2\sqrt{rr'}} 
\left( \Chi \sqrt{\frac{2}{1+\Chi}} 
\,  K \left( \sqrt{\frac{2}{1+\Chi}} \right) - 
\sqrt{2(1+\Chi)} 
\, E\left( \sqrt{\frac{2}{1+\Chi}} \right)
\right).
\end{equation}
Gradients of $Q_{\frac{1}{2}}$, and subsequently $g_1$,
 can be calculated analytically using the above
formulas as well.

\subsection{Discretization of integral operators}

Several geometries relevant to magnetic confinement in toroidal devices are
smooth, analytically parameterized, and can be described
using a modest number of equispaced discretization
points~\cite{cerfontaylor,hudson}. We therefore apply the Nystr\"om
scheme  and choose to discretize the integral equations~\eqref{eq_inteq0}
at equispaced points in the arc length parameter $s$.
Other parameterization variables $t$ can be used, and the subsequent
quadrature formula adjusted by a factor $ds/dt$.
We use $16^\text{th}$-order Hybrid Gauss-trapezoidal quadrature rules~\cite{alpert}
to evaluate the singular integrals. See~\cite{hao_2014} for a
nice review of Nystr\"om discretization options.

After applying an $n$-point Nystr\"om discretization to~\eqref{eq_inteq0}, 
we are left with a linear system:
\begin{equation}\label{eq_nystrom}
\frac{1}{2} \sigma_0^j + \left(\sum_{k = 1}^n w_{ij} \, G(s_j,s_k) \, 
\sigma_0^k \right) + \alpha \, u_j= 0, \qquad \text{for } j=1,\ldots,n,
\end{equation}
where $G$ is a linear combination of $g_0$ and $g_1$, $\sigma_0^j \approx
\sigma(s_j)$, $s_j = (j-1)L/n$, and $u_j$ is computed from the
contribution of the harmonic vector field $\bmm_H$ to the potential.
The flux constraint can be calculated as:
\begin{equation}
  \begin{aligned}
    \int_{S_t} \bB \cdot d\ba &= \frac{1}{\lambda} \int_\gamma \bB
    \cdot d\btau \\
&= \frac{1}{\lambda} \int_o^L \left( i\lambda \btau \cdot \cS_0 [\bmm_0] +
  i \btau \cdot \nabla \times \cS_0[\bmm_0] \right) ds
  \end{aligned}
\end{equation}
This integral can be discretized using the same Nystr\"om
scheme as in~\eqref{eq_nystrom}.

The above formula corresponds to the genus-one case; the genus-two
case follows exactly, with the poloidal flux being calculated as:
\begin{equation}
  \begin{aligned}
    \int_{S_p} \bB \cdot d\ba &= \frac{1}{\lambda} \int_{\partial S_p}
    \bB \cdot d\bl \\
&= \frac{1}{\lambda} \left( \int_0^{2\pi} B_0^\phi(s_\myout,\phi) \, d\phi 
- \int_0^{2\pi} B_0^\phi(s_\myin,\phi) \, d\phi  \right) \\
&= \frac{2\pi}{\lambda} \left(  B_0^\phi(s_\myout) - B_0^\phi(s_\myin)  \right),
  \end{aligned}
\end{equation}
where $B_0^\phi$ denotes the $\phi$-component of the purely
axisymmetric mode of the field $\bB$. The
integrals reduce to constants since $\bB_0$ is invariant as a function
of $\phi$, and $s_\myout$ and $s_\myin$ denote values along $\gamma$
which form $\partial S_p$.
The resulting  discretized block-system is given as:
\begin{equation}
\begin{bmatrix}
\mtx{A}   &   \vct{u}_1 &\vct{u}_2 \\
\vct{v}^T_1 &  c_{11} & c_{12} \\ 
\vct{v}^T_2  &  c_{21} & c_{22}
\end{bmatrix}
\begin{bmatrix}
\boldsymbol{\sigma} \\
\alpha \\
\beta
\end{bmatrix}
= \begin{bmatrix}
\boldsymbol{0} \\
\fluxtor \\
\fluxpol
\end{bmatrix}.
\end{equation}

\subsection{Axisymmetric surface differential operators}

In order to numerically apply the operator $\cK$, the 
surface gradient $\nabla_\Gamma$ and  the inverse surface Laplacian
$\surflap^{-1}$ need to be applied. 
Given an equispaced, periodic discretization in arc length 
of the generating curve $\gamma$, and a similar sampling of the function
$\sigma_0$, these operators can be
applied spectrally using Fourier methods in the arc length parameter
$s$.
See~\cite{nedelec,frankel_2011} for a general discussion of these
operators on arbitrary surfaces in three dimensions.

Recalling the discussion of the axisymmetric integral equation and
discretization  in the previous sections,
we only need to derive expressions for the surface
differential operators when applied to functions which have no angular
dependence. 
Let $f = f(s)$ be a scalar
function defined on $\Gamma$.
 Then, relative to the local
coordinate system $\btau$, $\bphi$, $\vct{n}$, we have the following
three expressions for the surface differentials:
\begin{equation}\label{eq-surfops}
\begin{gathered}
\nabla_\Gamma f = \frac{d f}{d s} \btau, \\
\surflap f = \frac{d^2 f}{d s^2} + 
\frac{1}{r} \frac{d r}{d s} \frac{d f}{d s}.
\end{gathered}
\end{equation}
If the generating curve were not parameterized in arc length, but
instead in some other variable $t$, the
previous formulae would include terms involving $ds/dt$ and $d^2s/dt^2$.

Fourier differentiation can easily be used to numerically compute the
derivatives $d/ds$. 
Due to the periodicity of $\gamma$, $f$ can be written as
\begin{equation}
f(s) = \sum_m f_m \, e^{2\pi i m s/L}.
\end{equation}
where $L$ is the length of $\gamma$. Let $\vct{f}$ be the discretization of the
function $f$ at equispaced points on $\gamma$ given by $s_1$, \ldots,
$s_n$.
The approximation $\vct{f}' \approx df/ds$ can then be computed as 
$\vct{f}' = \mtx{D} \vct{f}$, where $\mtx{D}$ is a Fourier spectral
differentiation matrix~\cite{trefethen_2000}. Linear combinations,
compositions, and row-scalings of $\mtx{D}$ can then be taken to
implement each of the above surface differentials. Let $\mtx{G}$ be
the matrix approximating $\surfgrad$ and $\mtx{L}$ be the matrix
approximating $\surflap$.

While $\mtx{L}$ can be directly
constructed using $\mtx{D}$, its inverse $\mtx{L}^{-1}$ requires some
more care. In fact, $\mtx{L}^{-1}$ does not exist since $\mtx{L}$ is
rank-one deficient (i.e. $\surflap c = 0$, where $c$ is a constant).
However, the surface Laplacian is uniquely invertible when restricted
as a map from mean-zero functions to mean-zero functions~\cite{EpGr}.
As shown in~\cite{sifuentes}, we can use this fact to instead solve
the integro-differential equation
\begin{equation}
\left( \frac{d^2 }{d s^2} + 
\frac{1}{r} \frac{d r}{d s} \frac{d }{d s}
+ \int_\Gamma  \right) \rho = f.
\end{equation}
This equation is invertible, and the solution $\rho$ 
has the properties that
$\surflap \rho = f$ and $\int_\Gamma \rho = 0$.

Analogously, we find $\mtx{L}^{-1}$ restricted to mean-zero functions
by inverting the matrix 
\begin{equation}
\mtx{L}_0 = \mtx{L} + \vct{1} \vct{w}^T,
\end{equation}
with $\vct{w}$ a vector of quadrature weights such that $\vct{w}^T \vct{u}
\approx \int_\gamma u$, and $\vct{1}$ is a vector of all
ones. As shown in~\cite{sifuentes}, $\mtx{L}_0^{-1}$ exists since
$\vct{1}$ is in the null-space of $\mtx{L}^T$ 
 (i.e. the null-space of the surface Laplacian is constant functions).

\subsection{Surface harmonic vector fields}

In the representation given in~\eqref{eq-mbuild}, 
it was assumed that the harmonic vector
field $\bmm_H$ was known a priori, as it was only the coefficient $\alpha$
that was unknown. Fortunately, along surfaces of revolution, one
can explicitly construct such harmonic vector
fields~\cite{werner_1966,EpGrOn}.
A basis for the two-dimensional
surface harmonic  vector fields along $\Gamma$, 
denoted by $\bmm_{H_1}$ and $\bmm_{H_2}$, is
\begin{equation}
\bmm_{H_1} = \frac{1}{r} \btau, \qquad
\bmm_{H_2} = -\frac{1}{r} \bphi,
\end{equation}
where we have used the relationship $\bmm_{H_2} = \vct{n} \times
\bmm_{H_1}$.
Note that the basis of harmonic vector fields on a genus-one surface is
two-dimensional.  For the purposes of constructing Beltrami fields in
genus-one toroidal volumes, we
must find a {\em single} surface harmonic vector field which 
satisfies \mbox{$\vct{n} \times \bmm_H =  i \bmm_H$}, see 
relation~\eqref{eq-surfharm}.
It is easy to check that the field
\mbox{$\bmm_H = \bmm_{H_1} -i \bmm_{H_2}$} satisfies this criterion.

In the genus-two case it is necessary to construct
harmonic vector fields along $\Gamma^\myout$ and $\Gamma^\myin$ such
that
\begin{equation}
\begin{aligned}
\vct{n} \times \bmm_H^\myout &= i  \bmm^\myout_{H}, \\ 
\vct{n} \times \bmm_H^\myin &= -i \, \bmm^\myin_{H},
\end{aligned}
\end{equation}
as described in Section~\ref{sec_shell}. In this case, it is the
fields \mbox{$\bmm_H = \bmm_{H_1} \mp i \bmm_{H_2}$}, respectively,
that satisfy these conditions. The above conditions are analogous to
those in~\eqref{eq-surfharm}, where the $\pm$-sign comes from the
definition of the normal vector~$\vct{n}$.

\subsection{A direct solver}

Due to the efficiency of having decomposed the surface integral
equation via Fourier methods into an integral equation along a
one-dimensional curve in the $rz$-plane, relatively high-order
accurate solutions can be obtained using a modest number of
unknowns. The resulting linear systems can be solved directly using
dense linear algebra in a few seconds. All of the code for the
following examples has been implemented in Fortran 90, and all linear
systems were solved densely using the Intel MKL implementation of
Lapack~\cite{lapack} LU-factorization routines.  The computational
cost of the solve is $\mathcal O(N^3)$, where $N$ is the number of
discretization nodes along the generating curve ($N$ is at most a few
hundred in our examples).  Timings are not reported, as the size of
these linear systems is extremely small (but yet are able to achieve
very high accuracy).  Specialized Fortran subroutines have been
written for every aspect of the code, including the modal Green's
function evaluation and special function evaluation. No third-party
software package were used, except Lapack.  For more complicated
geometries, or for very large values of the Beltrami parameter
$\lambda$, fast direct solvers in two-dimensions~\cite{young} or
non-translation invariant FMM-type methods would need to be used.

In the genus-one case, 
denoting by $\cS_0$ and $\cK_0$ the axisymmetric modal operators
 of the layer potential 
operators $\cS$ and $\cK$ from~\eqref{eq-inteq1}, 
the system matrix can be directly constructed as
\begin{equation}
\begin{aligned}
\begin{bmatrix}
\frac{1}{2}\mathcal I + \mathcal K_0 \\
\int_{\partial S_t} \bB \cdot d\bl
\end{bmatrix}
(\sigma,\alpha)
&\approx 
\begin{bmatrix}
\frac{1}{2} \mtx{I} + \mtx{S}' + i\lambda \mtx{N} \mtx{S} \mtx{M}
+ i \mtx{N} \mtx{S}^\times \mtx{M} & \vct{u} \\
\vct{w}^T(i\lambda \mtx{N} \mtx{S} \mtx{M}
+ i \mtx{N} \mtx{S}^\times \mtx{M}) & c
\end{bmatrix}
\begin{bmatrix}
\boldsymbol{\sigma} \\
\alpha
\end{bmatrix}= \begin{bmatrix}
\boldsymbol{0} \\
\Phi_{t}
\end{bmatrix}\\
\end{aligned}
\label{eq:large_system}
\end{equation}
where
\begin{equation}
\begin{aligned}
\mtx{S} &\approx \mathcal S_0, &\qquad \mtx{N} &\approx \vct{n} \cdot, \\
\mtx{S}' &\approx \mathcal S'^*_0, & 
   \mtx{M} &= i\lambda\left( \mtx{G}\mtx{L}_0^{-1} - i \mtx{N}^\times \mtx{G}
\mtx{L}_0^{-1} \right), \\
\mtx{S}^\times &\approx \nabla \times \mathcal S_0, &  
\mtx{N}^\times &\approx \vct{n} \times,
\end{aligned}
\end{equation}
and where $c$ is the contribution to the toroidal flux due to
$\bmm_{H}$. Similarly, we can construct the corresponding system
matrix for the genus-two case using the same set of matrices as above,
with the additional poloidal flux constraint.

Detailed comparisons between the performance of our new
solver and existing solvers, such as the Beltrami solver in the code
SPEC~\cite{hudson}, would be valuable and are planned for the near
future. At this stage, we can mention two desirable features of our
solver which are independent of the details of the actual implementation of
the numerical scheme, and follow directly from the particular
integral formulation in
which only the boundary of the domain needs to be discretized. First,
we are not faced with artificially-induced
issues associated with the potential existence
of a coordinate singularity which naturally occur when discretizing
the volume of genus-one domains~\cite{hudsonPPCF}. Second, the number
of unknowns in our solver is much smaller, which leads to significant
savings in memory usage. As an illustration, using a standard of
$N_{s}=16$ Chebyshev nodes to discretize the volume in the radial
direction, the number of unknowns in our scheme is 16 times smaller 
than it is in a volume-based solver. These savings in memory usage are
even more dramatic in situations in which a global magnetohydrodynamic
(MHD) equilibrium is computed by subdividing the domain into a large
number of genus-one and genus-two toroidal regions which are each in a
Taylor state~\cite{hudson}. In such situations, our solver would solve
for $N_{s}\times N_{V}$ fewer unknowns than a volume based solver
would, where $N_{V}$ is the number of toroidal Taylor state regions
used to describe the MHD equilibrium. A standard resolved calculation
may use $N_{s}=16$ and $N_{V}=10$~\cite{hudson}, so that our scheme
would have 160 times fewer unknowns.

On the other hand, it is harder to predict the
advantages of our scheme regarding run times. This can be explained as
follows. In certain formulations (e.g. volume finite-difference or
finite-element schemes),
the matrices for the linear system
which results from discretizing $\nabla\times\bB = \lambda\bB$ are
sparse~\cite{hudson}, and inverting the resulting system is not
significantly slower than inverting our smaller, dense
system~\eqref{eq:large_system}. In
volume based solvers, $\bB$ is solved for at every grid point in the
volume, whereas in our formulation, the field still needs to be
evaluated at the desired points in the volume using numerical
quadratures for evaluating the expressions in~\eqref{eq-pot_convol}. This
suggests that, as far as asymptotic run-times are concerned,
our approach is most
promising in situations in which the magnetic field only needs to be
known on the boundary of the domain \cite{hudson}. In any case, this
discussion motivates detailed code performance comparisons with SPEC,
which we are planning to do in the near future.

Before closing this section, we observe that the numerical
method to compute Taylor states we present in this article is not
restricted to domains with a smooth boundary. The generalized
Debye representation used in our solver is also applicable to
regions with corners~\cite{cherno}, and therefore
our scheme is also expected to be
effective for situations in which the boundary has one or several
corners. These types of geometries correspond
to the existence of one or several magnetic
X-points~\cite{cerfontaylor}. There are, however, slight technical
changes that would need to be made.
The first complication is that for domains with a
corner, the discretization of the integral operators would have to be
done using a panel-based discretization scheme, as used by 
Bremer in~\cite{bremer_2010,bremer_2012} instead of our hybrid
Gauss-trapezoidal rule. Likewise, Fourier based differentiation would
have to be replaced with panel based spectral differentiation. While
these schemes are available, we have not yet implemented them, and are
therefore not able to report on the effectiveness of our scheme for
these situations.

\section{Examples}
\label{sec_examples}

In this section we consider the following three numerical examples,
which we use to evaluate the accuracy of our solver:
\begin{enumerate} 
\item A shaped, low aspect ratio Taylor state as may be
observed in spheromak experiments~\cite{geddes},
\item A Taylor state in a toroidal shell taken from a 
sub-region of the equilibrium in the previous example, and
\item A sequence of non-axisymmetric Taylor states in a shaped,
  axisymmetric, moderate aspect ratio domain $\Omega$.
\end{enumerate}

\subsection{Example 1: Axisymmetric Taylor state of a 
shaped plasma with very low aspect ratio}
In a recent article~\cite{cerfontaylor}, we described a general method
to compute exact analytic Taylor states in shaped, axisymmetric
plasmas. The method of that work 
is used in this section to evaluate the point-wise
numerical error obtained by our solver.  We start with a
quick review of the construction of the analytic Taylor state before
presenting our numerical results. The construction relies on writing
the magnetic Beltrami field in $\Omega$ that satisfies
$\nabla\times\boldsymbol{B}=\lambda\boldsymbol{B}$ as
\begin{equation}
\bB = \lambda\frac{\psi}{r} \bphi
  + \frac{1}{r} \nabla\psi\times \bphi,
\label{Bexact}
\end{equation}
where~$\psi$ is the poloidal flux function. This function
$\psi$ satisfies the Grad-Shafranov equation
\begin{equation}  \label{GSpsi}
  \begin{aligned}
      r\frac{\partial}{\partial
  r} \left(\frac{1}{r} \frac{\partial\psi }{\partial
  r}\right) + \frac{\partial^{2}\psi}{\partial z^{2}}
     &=-\lambda^{2}\psi & \qquad &\text{in } \Omega,\\
    \psi&=0 & &\text{on } \Gamma.
  \end{aligned}
\end{equation}
The boundary $\Gamma$ is undetermined for the moment.
Consider a general solution to~\eqref{GSpsi} given by
\begin{equation}\label{eq:general_sol}
\psi(r,z) = \psi_{0}
+ \sum_{i = 1}^5 c_{i}\, \psi_{i}
\end{equation}
with $c_i$ constants, and where the functions $\psi_{i}$ are defined by
\begin{align*}
\psi_{0}(r,z)&=rJ_{1}(\lambda r),\\
\psi_{1}(r,z)&=rY_{1}(\lambda r),\\
\psi_{2}(r,z)&=rJ_{1}\left(\sqrt{\lambda^{2}-c_{6}^{2}}r\right)\cos(c_{6}z),\\
\psi_{3}(r,z)&=rY_{1}\left(\sqrt{\lambda^{2}-c_{6}^{2}}r\right)\cos(c_{6}z),\\
\psi_{4}(r,z)&=\cos\left(\lambda\sqrt{r^{2}+z^{2}}\right),\\
\psi_{5}(r,z)&=\cos\left(\lambda z\right).
\end{align*}
with $J$ and $Y$ the Bessel functions of the first and second kinds,
respectively.
\begin{figure}[!t]
  \centering
  \begin{subfigure}[b]{.3\linewidth}
    \centering
    \includegraphics[width=.95\linewidth]{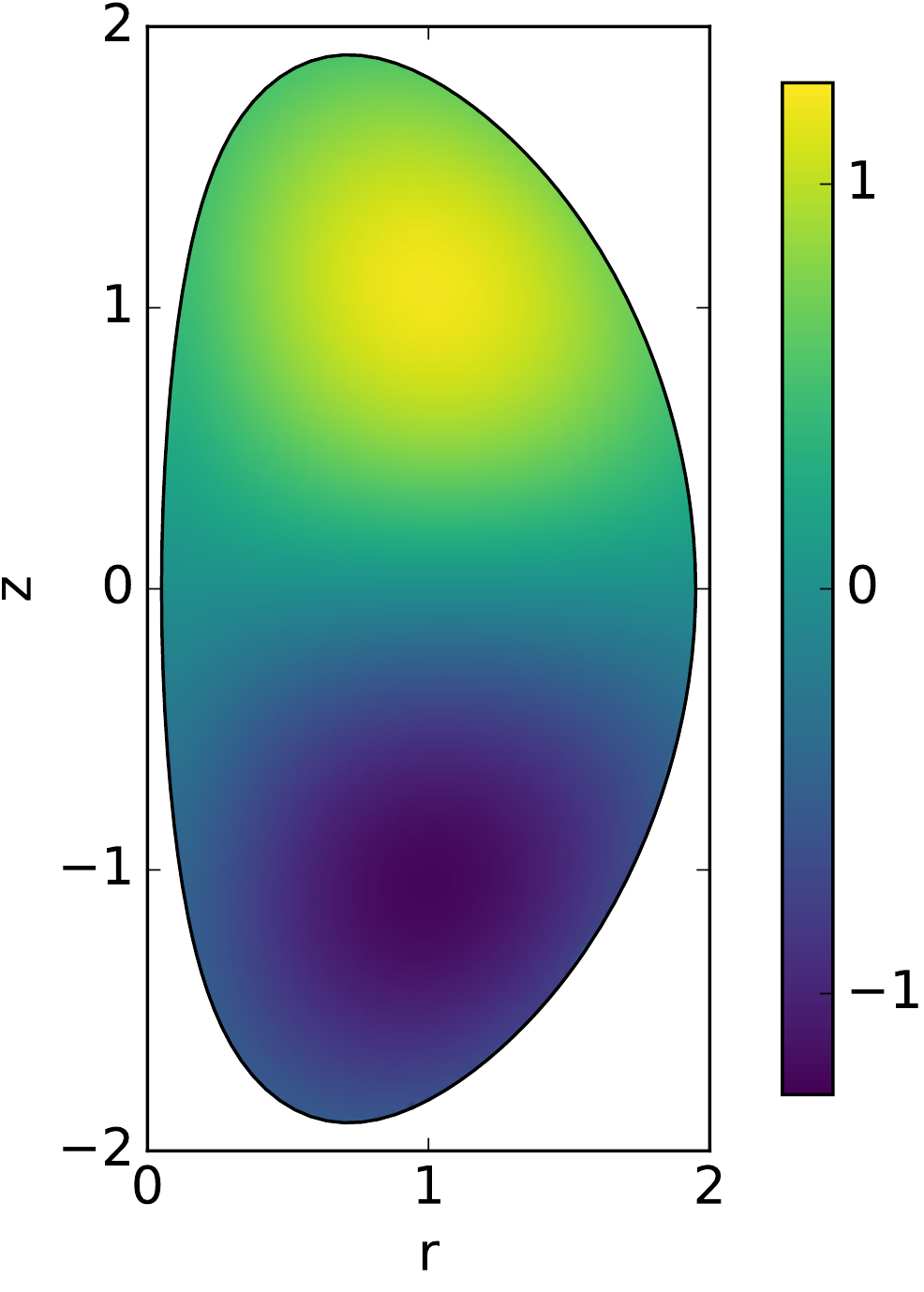}
    \caption{The radial component, $B_r$.}
    \label{fig_ex1r}
  \end{subfigure} \quad
  \begin{subfigure}[b]{.3\linewidth}
    \centering
    \includegraphics[width=.95\linewidth]{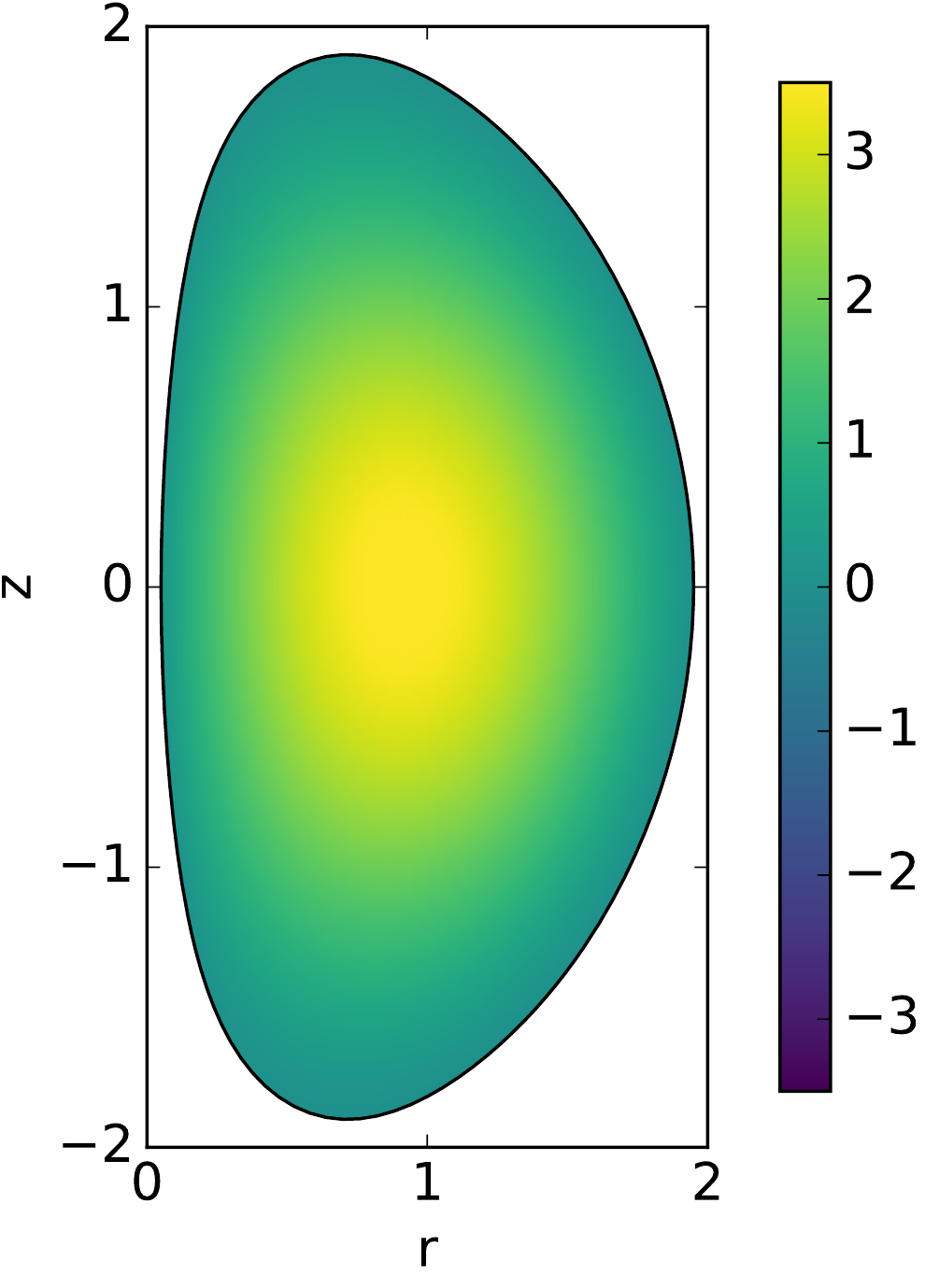}
    \caption{The azimuthal component, $B_\phi$.}
    \label{fig_ex1t}
  \end{subfigure} \quad
  \begin{subfigure}[b]{.3\linewidth}
    \centering
    \includegraphics[width=.95\linewidth]{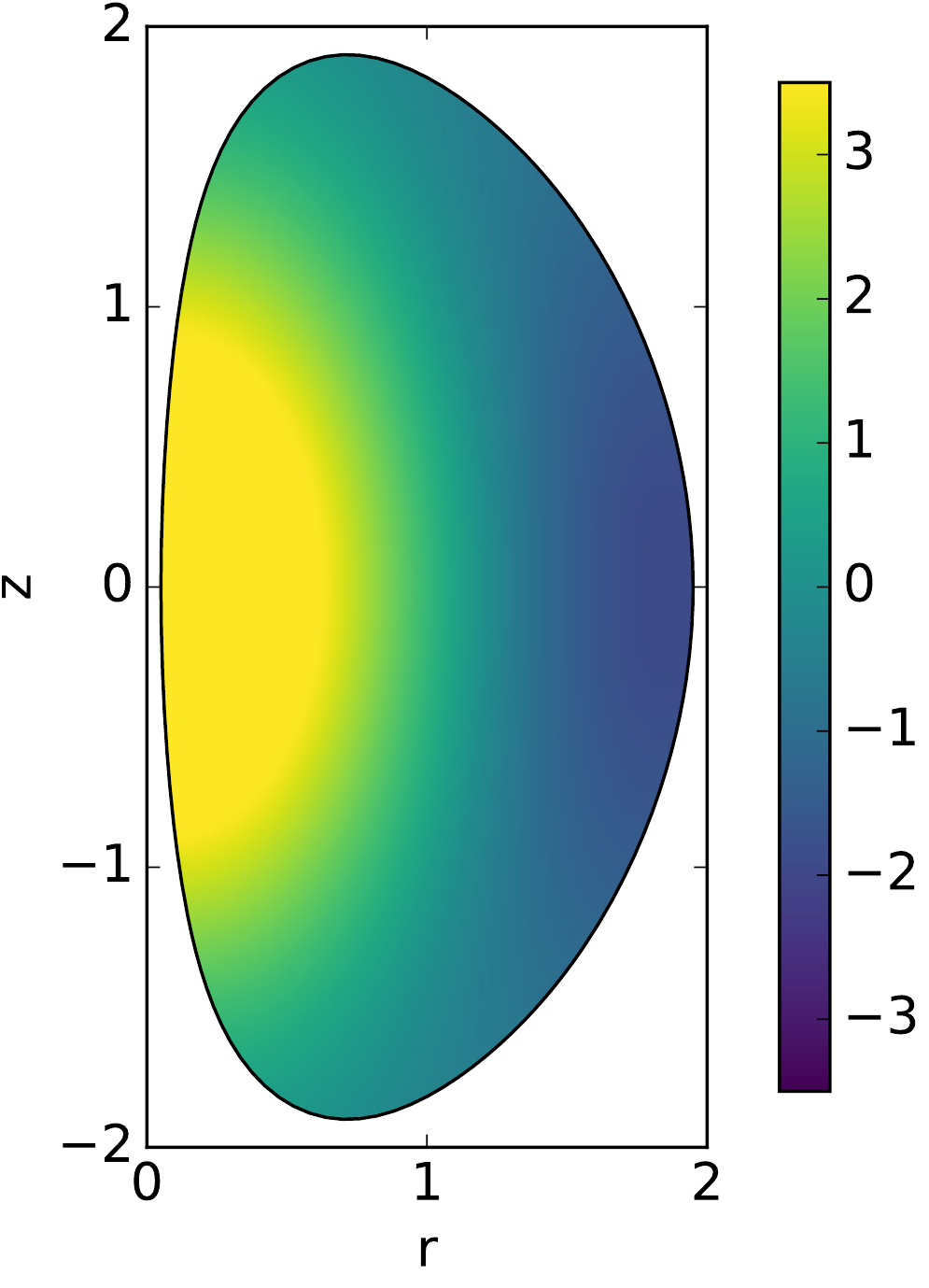}
    \caption{The $z$ component, $B_z$.}
    \label{fig_ex1z}
  \end{subfigure}
  \caption{The Beltrami field $\bB$ in the $\phi=0$
    plane.}
  \label{fig_genus1}
\end{figure}
Treating~$\lambda$ and $c_1, \ldots, c_6$ as unknowns, we solve
for these values by imposing seven boundary conditions on $\psi$
chosen such that the boundary of the plasma, given by the implicit
equation $\psi=0$, best approximates the parametric curve
\begin{equation}\label{eq:boundary}
  \begin{aligned}
    r(t)&=1+\epsilon  \cos\left(t+\alpha  \sin t \right)\\
    z(t)&=\epsilon \kappa  \sin t,
  \end{aligned}
\end{equation}
with $t \in [0,2\pi)$.
This curve describes a wide class of experimentally relevant axisymmetric
plasma boundaries \cite{miller}, where 
 $\epsilon$ is the inverse aspect ratio, 
$\kappa$ is the elongation, and $\sin \alpha=\delta$ is the
triangularity. 
The boundary conditions enforced on $\psi$ are~\cite{cerfonsolovev}
\begin{equation}\label{eq:constraints}
\begin{aligned}
\psi(1+\epsilon,0,\bc)&=0,\\
\psi(1-\epsilon,0,\bc)&=0,\\
\psi(1-\delta\epsilon,-\kappa\epsilon,\bc)&=0,\\
\psi_{r}(1-\delta\epsilon,-\kappa\epsilon,\bc)&=0,\\
\psi_{zz}(1+\epsilon,0,\bc)+N_{1}\psi_{r}(1+\epsilon,0,\bc)&=0,\\
\psi_{zz}(1-\epsilon,0,\bc)+N_{2}\psi_{r}(1-\epsilon,0,\bc)&=0,\\
\psi_{rr}(1-\delta\epsilon,-\kappa\epsilon,\bc)
+ N_{3}\psi_{z}(1-\delta\epsilon,-\kappa\epsilon,\bc)&=0,
\end{aligned}
\end{equation}
where the subscripts refer to partial derivatives with respect 
to the specified variable, 
the vector of parameters is $\bc=(c_{1}, \cdots, c_{6}, \lambda)$, 
and $N_{1}$, $N_{2}$ and $N_{3}$ are the curvatures at the outer 
equatorial point, the inner equatorial point, and the top 
point, respectively, given by the formulae
\begin{equation}
N_{1}=-\frac{(1+\alpha)^2}{\epsilon\kappa^{2}}, \qquad 
N_{2}=\frac{(1-\alpha)^2}{\epsilon\kappa^{2}}, \qquad 
N_{3}=\frac{\kappa}{\epsilon\cos^{2}\alpha}.
\end{equation}
Equation \eqref{eq:constraints} is a non-linear system of 7 equations
for 7 unknowns. Given reasonable initial conditions, it can be solved
without difficulty using standard non-linear root finding
packages.
\begin{table}[!b]
  \centering
  \caption{Convergence data for Example 1: A genus one Taylor
    state. Relative errors are reported in the field $\bB$ 
    at the point $r=1.2$,  $\phi=0$, $z=0.25$. The reference field
    $\bB_{\text{exact}}$ is calculated using $\psi$ from~\eqref{Bexact}.}
  \label{tab_genus1}
  \begin{tabular}{|c|ccc|c|}
    \hline
    $n$ & $\Re B_r$ & $\Re B_\phi$ & $\Re B_z$ & 
       $|\bB - \bB_{\text{exact}}|/|\bB_{\text{exact}}|$ \\
    \hline
    25 & 0.443052524078644 & 3.10056763474524 &
              -3.784408049008867E-002 & $2.7 \cdot 10^{-3}$ \\
    \hline
    50 & 0.442014263551259 & 3.09845144534915 &
              -4.109405171821609E-002 & $2.5 \cdot 10^{-5}$  \\
    \hline
    100 & 0.442018001760211 & 3.09850436011175 & 
             -4.104126312770094E-002  & $3.9 \cdot 10^{-8}$ \\
    \hline
    200 & 0.442017994270342 & 3.09850428092008 &
             -4.104130814605825E-002 & $1.2 \cdot 10^{-8}$ \\ \hline
  \end{tabular}
\end{table}
We achieve near machine precision residual using
Matlab's \texttt{fsolve} function.
Once $\bc$ has been determined,  the Taylor state is
completely defined. Using the fact that $\nabla \times \bB = \lambda
\bB$, and relation~\eqref{Bexact},
we can numerically evaluate the toroidal flux as
\begin{equation}\label{exactflux}
  \begin{aligned}
    \fluxtor &= \int_{S_t} \bB \cdot d\ba \\
    &= \frac{1}{\lambda} \int_\gamma \bB \cdot d\bl \\
    &= \frac{1}{\lambda} \int_\gamma \frac{1}{r} \left( \nabla \psi
      \times
    \bphi \right) \cdot d\bl.
  \end{aligned}
\end{equation}
This integral can be computed to spectral accuracy using the
trapezoidal rule.
The integral equation-based Beltrami solver of this paper 
can be tested against this exact Taylor state by solving the problem
\begin{equation}
  \begin{aligned}
    \nabla \times \bB &= \lambda \bB & \quad &\text{in } \Omega, \\
    \bB \cdot \vct{n} &= 0 & &\text{on } \Gamma,\\
    \int_{S_t} \bB \cdot d\ba &= \fluxtor. & & 
  \end{aligned}
\end{equation}

To highlight the robustness and accuracy of our solver, we choose
equilibrium parameters that are particularly challenging from a
numerical point of view: very low aspect ratio, $\epsilon=0.95$,
moderate triangularity, $\delta=0.3$, and high elongation, $\kappa=2$.
See Figure~\ref{fig_genus1} for a plot of the real part of $\bB$, and
Table~\ref{tab_genus1} for convergence data. We report the relative
$\ell^2$-error at the arbitrary point $(r,\phi,z)=(1.2,0,0.25)$ for various values of $n$,
the number of points used in discretizing the geometry. Observe that in our method, the interior of $\Omega$ is not discretized, and we are free to evaluate the field at any point of the domain. We chose $(1.2,0,0.25)$ because the error at this point is representative of the error at any other point inside $\Omega$. The solver converges rapidly to a precision of roughly~$10^{-8}$.
This precision is most likely due to
a combination of the accuracy to which we evaluated
the modal-Green's functions -- set to be an absolute precision of
$10^{-12}$ inside the adaptive Gaussian quadrature routine -- and the
conditioning of the system matrix -- $\kappa(\mtx{A}) \sim 10^4$
incurred by spectral differentiation.
Increased precision could be obtained by more accurate evaluation of
the modal Green's functions, and constructing the approximation to
$\surflap^{-1}$ using an integral method instead of inverting a
spectral differentiation matrix. These minor numerical improvements are
discussed more thoroughly in~\cite{epstein_2017}.

\subsection{Example 2: Axisymmetric Taylor state within a toroidal
  shell with very low aspect ratio}

The previous equilibrium can also be used to verify the accuracy and
performance of our solver for Taylor states in toroidal shells. The
flux surface $\psi=0.5$ can be taken as the inside surface,~$\Gamma^\myin$,
 and the surface $\psi=0$ can be taken as the outside 
surface,~$\Gamma^\myout$. The toroidal shell is then bounded by $\Gamma^\myin$
and $\Gamma^\myout$. Let $S_t$ be the spanning 
surface obtained from the intersection of $\Omega$ with the plane $\phi=0$.
In this case, the toroidal  flux is given as
\begin{equation}
    \fluxtor = \frac{1}{\lambda} \left( \int_{\gamma_\myout}
      \frac{1}{r} \left( \nabla \psi \times \bphi \right) \cdot d\bl -
      \int_{\gamma_\myin} \frac{1}{r} \left( \nabla \psi \times \bphi
      \right) \cdot d\bl \right),
\end{equation}
where $\gamma_\myout$ and $\gamma_\myin$ are the generating curves for
$\Gamma^\myout$ and $\Gamma^\myin$, respectively. Furthermore, by
definition of the flux function $\psi$~\cite{cerfontaylor}, 
the poloidal flux is given by the simple formula:
\begin{equation}
  \begin{aligned}
    \fluxpol &= 2\pi \left( \psi(r_\myout,0) - \psi(r_\myin,0) \right) \\
    &= 2\pi \left( 0 - 0.5 \right) \\
    &=-\pi
\end{aligned}
\end{equation}
The values $r_\myout$ and $r_\myin$ are the radial coordinates
corresponding to $\psi=0$ and $\psi=0.5$ for $z=0$, respectively. With
these conditions, the magnetic field $\bB$ is given by~\eqref{Bexact},
with $\psi$ given by~\eqref{eq:general_sol}, and the constants
$c_{1},\ldots,c_6$, and $\lambda$ determined by solving the
system~\eqref{eq:constraints}.

\begin{figure}[t!]
  \centering
  \begin{subfigure}[b]{.3\linewidth}
    \centering
    \includegraphics[width=.95\linewidth]{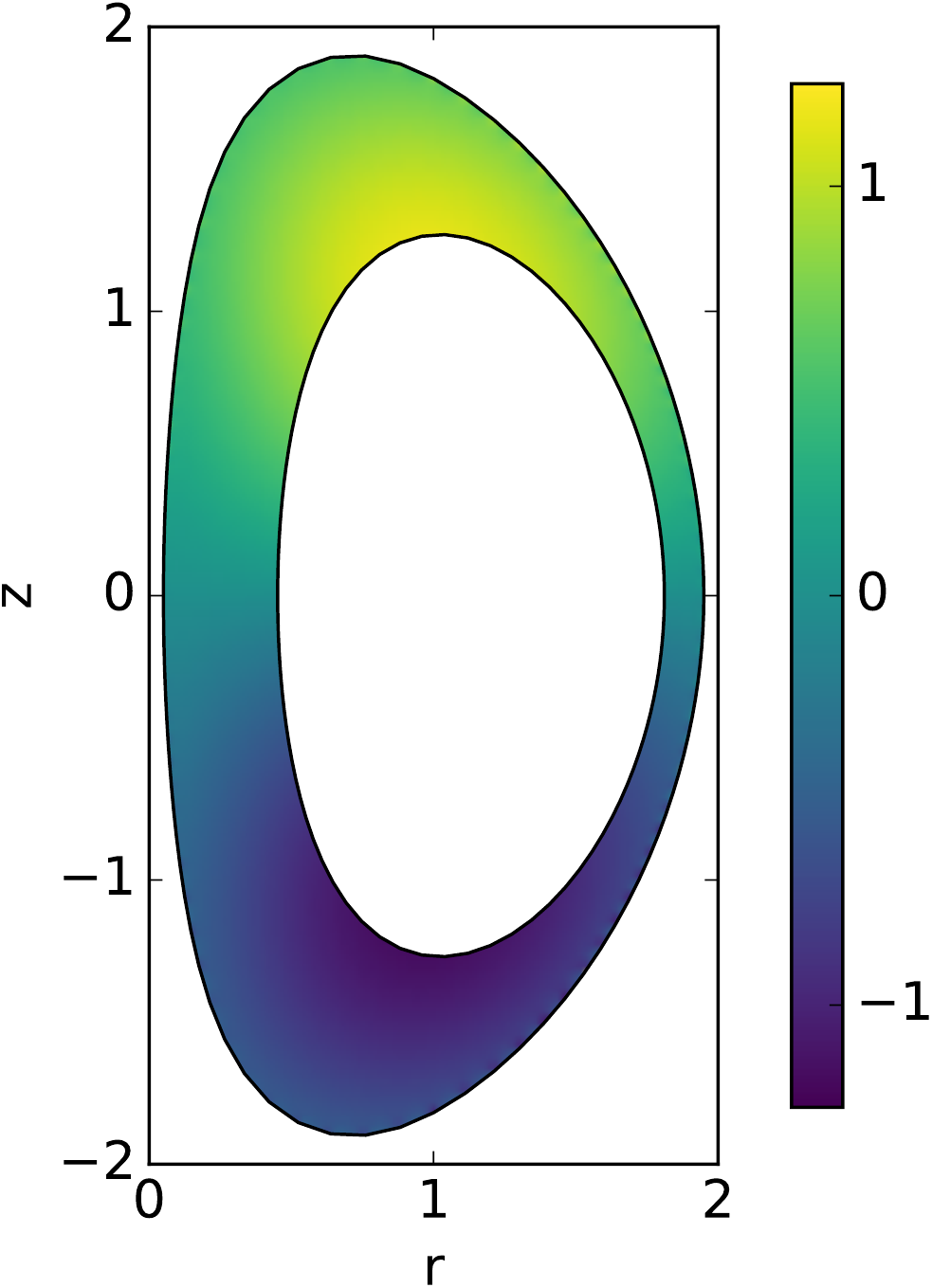}
    \caption{The radial component, $B_r$.}
    \label{fig_ex2r}
  \end{subfigure} \quad
  \begin{subfigure}[b]{.3\linewidth}
    \centering
    \includegraphics[width=.95\linewidth]{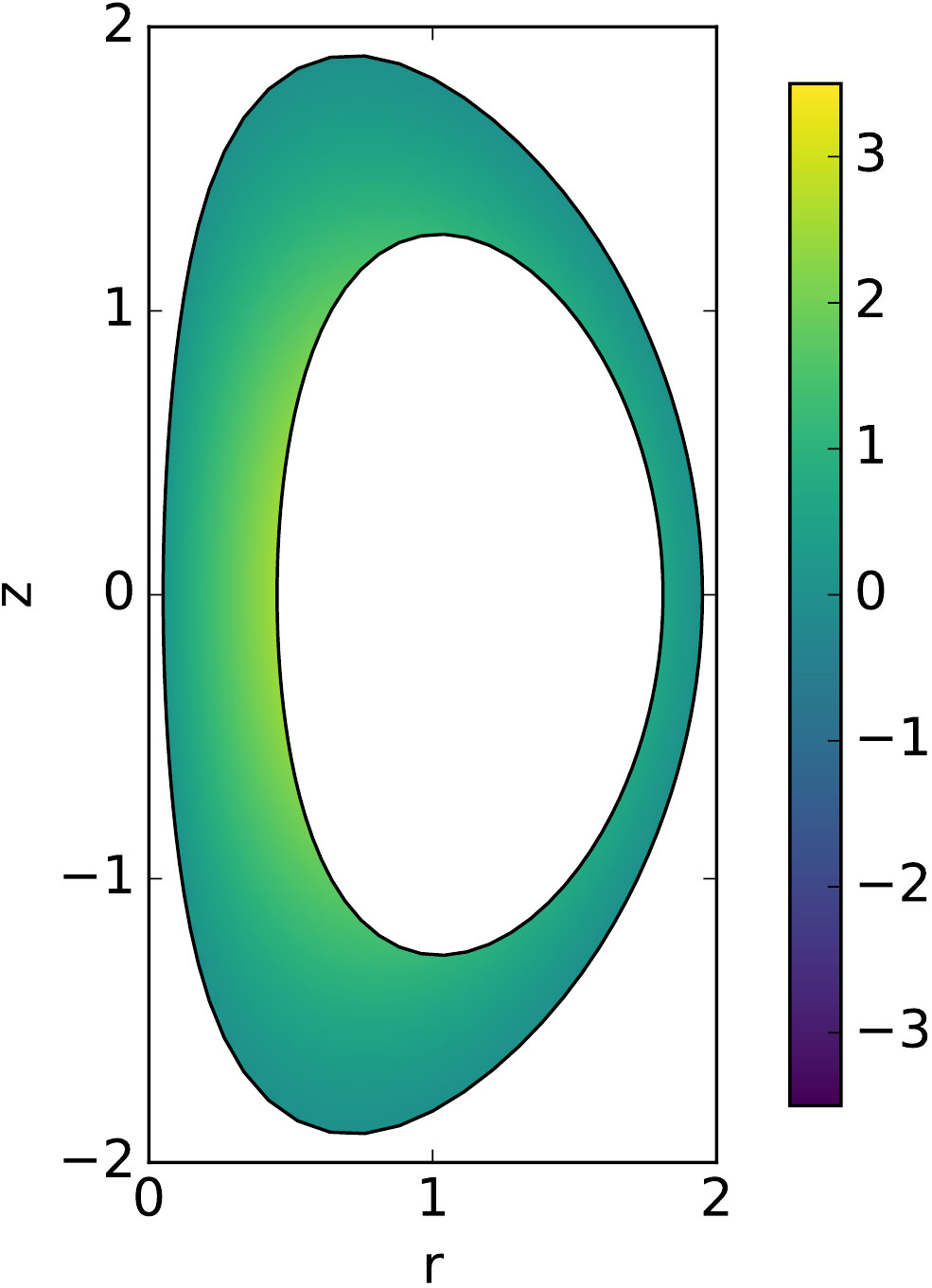}
    \caption{The azimuthal component, $B_\phi$.}
    \label{fig_ex2t}
  \end{subfigure} \quad
  \begin{subfigure}[b]{.3\linewidth}
    \centering
    \includegraphics[width=.95\linewidth]{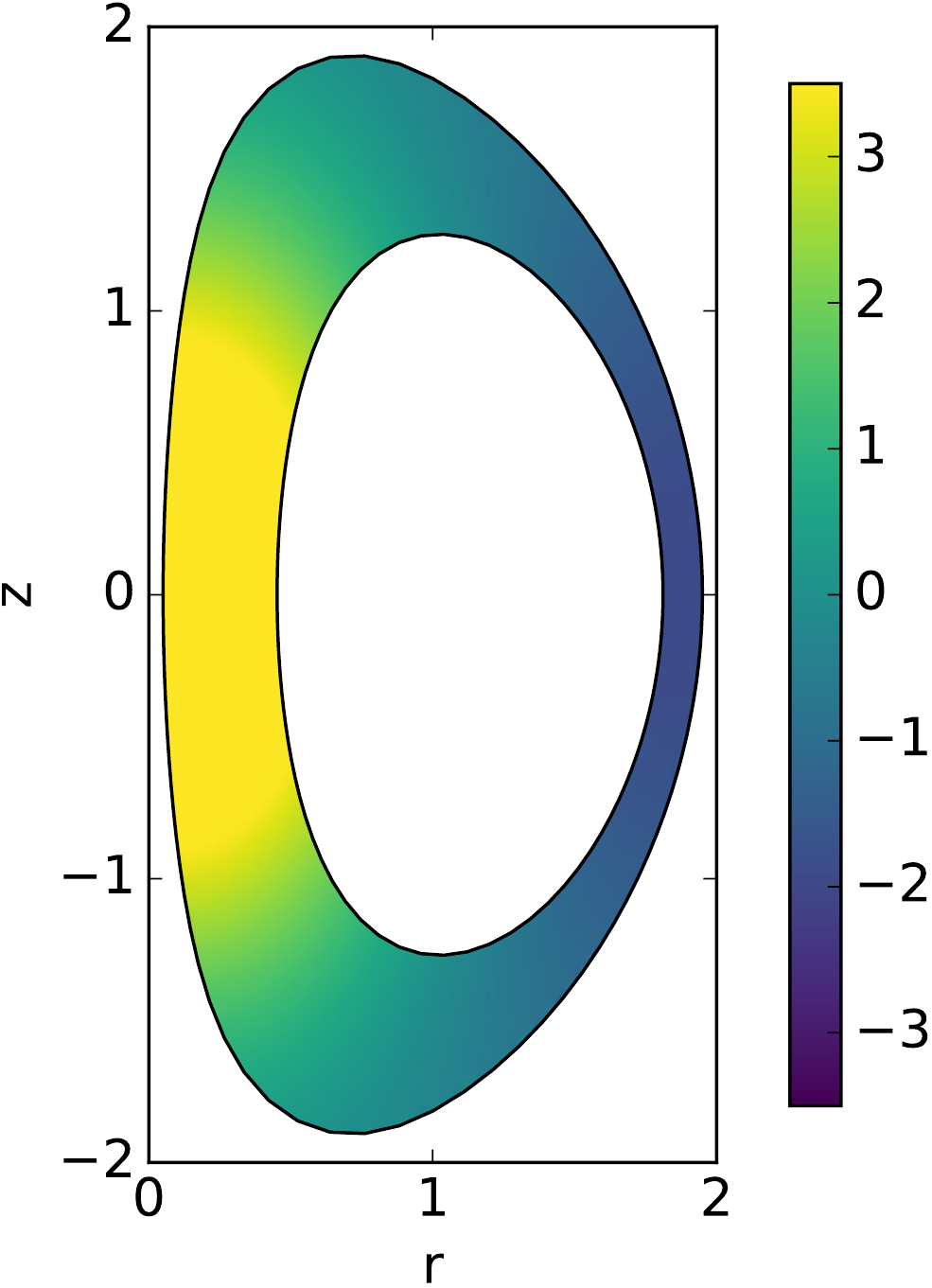}
    \caption{The $z$ component, $B_z$.}
    \label{fig_ex2z}
  \end{subfigure}
  \caption{The Beltrami field $\bB$ in the $\phi=0$ plane bounded by
    the surfaces $\psi=0.5$ and $\psi=0$.}
  \label{fig_genus2}
\end{figure}

See Figure~\ref{fig_genus2} for a plot of the real part of $\bB$ in
the toroidal shell, and
Table~\ref{tab_genus2} for convergence data. It is worth noting that
the field $\bB$ in the genus two shell 
is \emph{exactly} the same $\bB$ as in the genus one torus. The
additional boundary was merely calculated from the flux
function~$\psi$. 
However, the integral equation solver has no knowledge of
this, and is able to reproduce exactly the same field.
As before, we report the relative
$\ell^2$-error at the point $(r,\phi,z)=(1.2,0,0.25)$ for various values of $n$,
the number of points used in discretizing the geometry. We observe once more that our solver converges to a very accurate solution with few points on the boundary.

\begin{table}[!b]
  \centering
  \caption{Convergence data for Example 2: A genus two Taylor
    state. Relative errors are reported in the field $\bB$ 
    at the point $r=0.5$,  $\phi=0$, $z=-1.5$. The reference field
    $\bB_{\text{exact}}$ is calculated using $\psi$ from~\eqref{Bexact}.}
  \label{tab_genus2}
  \begin{tabular}{|c|ccc|c|}
    \hline
    $n$ & $\Re B_r$ & $\Re B_\phi$ & $\Re B_z$ & 
       $|\bB - \bB_{\text{exact}}|/|\bB_{\text{exact}}|$ \\
    \hline
    25  &-0.7790590773628590 & 0.5058371725845370 &
           0.9957374643442100 & $1.3 \cdot 10^{-2}$ \\ \hline
    50  &-0.7758504363890280 & 0.5043487336557070 &
           0.9869834030024680 &$7.3 \cdot 10^{-4}$ \\ \hline
    100 &-0.7754614741802320 & 0.5046765566326400 &
           0.9867586238268730 &$3.0 \cdot 10^{-6}$ \\ \hline
    200 &-0.7754611961412940 & 0.5046760189196530 &
        0.9867575110491060 &$8.6 \cdot 10^{-7}$ \\ \hline
  \end{tabular}
\end{table}

\subsection{Example 3: Non-axisymmetric Taylor state}

Solutions to the boundary value problem~\eqref{eq-bv} 
are not always purely axisymmetric. Only purely axisymmetric solutions
contribute to the toroidal and poloidal fluxes, but non-trivial, zero-flux Beltrami fields also exist~\cite{taylor_review}.
We focus here on the genus-1 toroidal case for simplicity, 
as the genus-2 case is identical.
Let $\lambda$ be an eigenvalue for the boundary value
problem corresponding to the $\lth$-Fourier mode solution. That
is to say, let~$\lambda$ be such that 
there exists a non-trivial solution to the problem:
\begin{equation}
\begin{aligned}
\nabla \times \bB_\ell &= \lambda \bB_\ell & \quad & \text{in } \Omega, \\
\bB_\ell \cdot \vct{n} &= 0 & & \text{on } \Gamma,
\end{aligned}
\end{equation}
where $\bB_\ell = \bB_\ell(r,z) \, e^{i\ell\phi}$.

Using the method of~\cite{zhao_2015,boyd_2002}, roots of the Fredholm determinant of the
Beltrami integral equation~\eqref{eq-inteq1} can be found as a function of the
Beltrami parameter~$\lambda$.  These roots correspond exactly to eigenvalues of the Beltrami boundary-value problem. The eigenvalues for the $\ell=1$ mode located on the interval
$[1.0,8.0]$ are reported in Table~\ref{tab_eig}, along with absolute
errors estimated as in~\cite{zhao_2015}. These values were obtained
from discretizing the $\ell=1$ Fourier mode, in a way analogous to the azimuthal
mode $\ell=0$ as described in Section~\ref{sec_numerical}, using $n=100$
points on the boundary of the shape given by the following parameterization
\begin{equation}
\begin{aligned}
 r(t) & =2 + \epsilon  \cos\left(t+\alpha  \sin t \right)\\
 z(t) & =\epsilon \kappa  \sin t,
\end{aligned}
\end{equation}
for $t \in [0,2\pi)$ with $\epsilon = 0.85$, $\sin\alpha = 0.3$, 
and $\kappa  = 2.0$.
The errors reported in Table~\ref{tab_eig} are generally commensurate
with the smallest singular value of the discretized system.

\begin{table}[!b]
  \centering
  \caption{Beltrami eigenvalues for the $\ell=1$ Fourier mode on the
    interval~$[1,8]$.}
  \label{tab_eig}
  \begin{tabular}{|c|c|}
    \hline
    $\lambda$ & Error \\ \hline
    2.81618429764383   &   $2.3\cdot 10^{-7}$    \\ \hline
     3.22821787079846  &   $  3.3\cdot 10^{-7}$   \\ \hline
     4.01342328856135  &   $ 5.8\cdot 10^{-9}  $ \\ \hline
     4.45732687692555   &  $  1.0\cdot 10^{-8} $  \\ \hline
     4.75909602398894  &   $  1.3\cdot 10^{-8} $  \\ \hline
     4.80160935115718  &   $  8.9\cdot 10^{-9} $  \\ \hline
     5.52819229381708  &   $ 1.0\cdot 10^{-10 }$  \\ \hline
     5.56546068190407  &   $ 5.8\cdot 10^{-10} $  \\ \hline
     6.13551340937516  &   $  6.4\cdot 10^{-12} $  \\ \hline
     6.34490415618171  &   $ 1.4\cdot 10^{-11} $  \\ \hline
     6.55792492108800  &   $  8.2\cdot 10^{-12}$ \\ \hline
     6.63664744243683  &   $  1.1\cdot 10^{-11}$ \\ \hline
     7.07387937977634  &   $  4.8\cdot 10^{-12}$ \\ \hline
     7.14679867372582  &   $  5.0\cdot 10^{-12}$ \\ \hline
     7.44941373173176  &   $  1.0\cdot 10^{-11}$ \\ \hline
     7.81008353287565  &   $  6.6\cdot 10^{-12}$ \\ \hline
     7.88508920256358  &   $  2.9\cdot 10^{-11}$ \\ \hline
  \end{tabular}
\end{table}

\begin{figure}[!t]
  \centering
  \includegraphics[width=.95\linewidth]{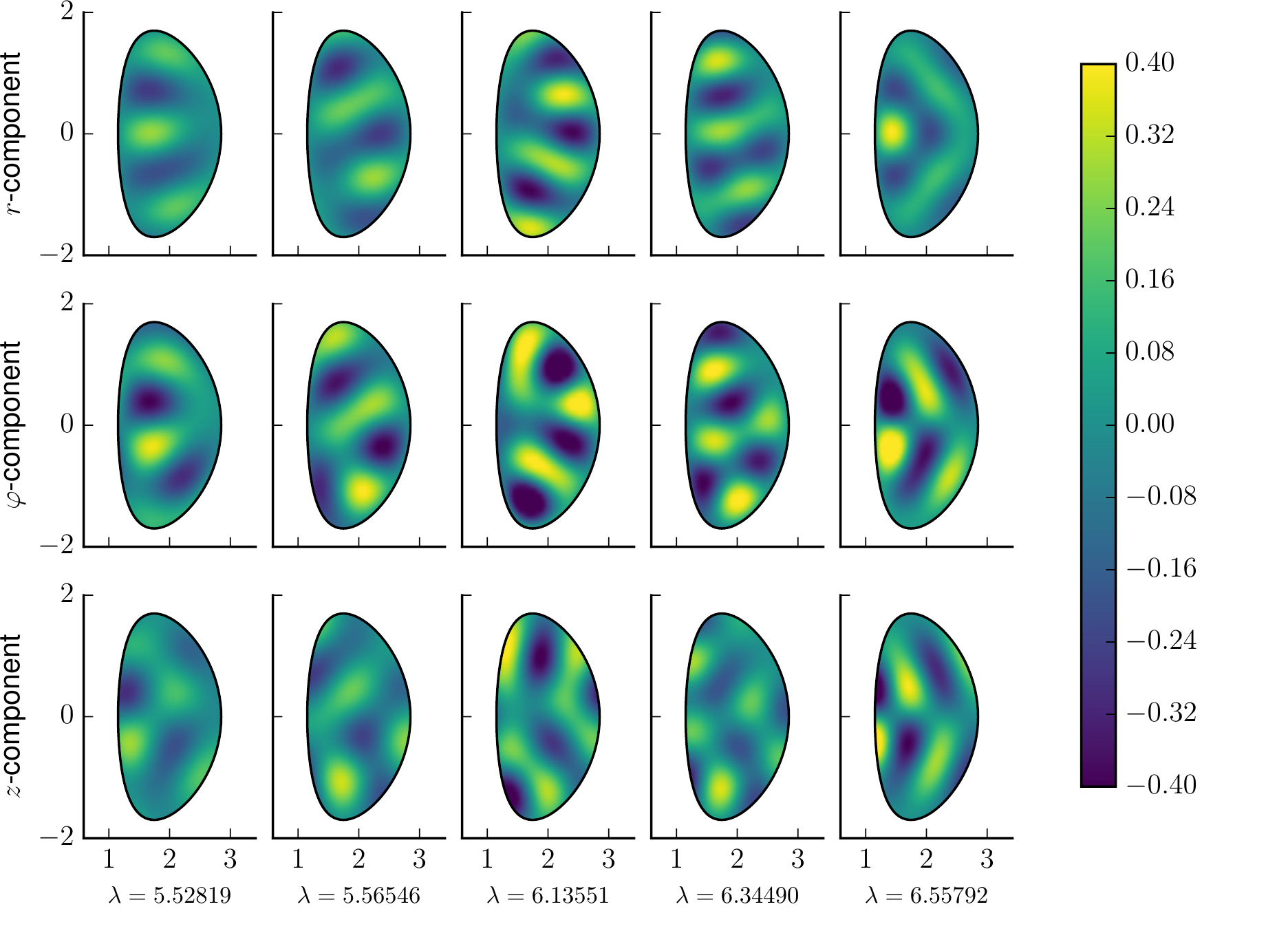}
  \caption{Slices of non-axisymmetric, zero-flux Beltrami fields
    in the $\phi=0$
    plane corresponding to the first five resonances larger than
    5 with azimuthal dependence of $e^{i\phi}$. 
    The modal Debye sources were normalized to have unit norm along
    the generating curve, 
    $\int_\gamma |\sigma_1|^2 ds= 1$.}
  \label{fig_eigs}
\end{figure}

Empirically, as shown via a numerical
calculation of the singular values of the discretized matrix, 
the eigenspace corresponding to each of the eigenvalues in 
Table~\ref{tab_eig} is one-dimensional. The eigenvector (i.e. the 
corresponding generalized Debye source) can be computed by using the method
of~\cite{sifuentes}. That is, we take a random vector
$\vct{r}$ and apply the discretized system matrix~$\mtx{A}$ to it to compute
$\vct{y} = \mtx{A} \vct{r}.$
Next, using standard Gaussian elimination, we solve the system
$\left( \mtx{A} + \vct{r}_1 \vct{r}^T_2\right) \vct{x} = \vct{y}$, with
    $\vct{r}_1$ and $\vct{r}_2$ random vectors. The matrix $\mtx{A} + 
\vct{r}_1 \vct{r}^T_2$ is uniquely invertible with probability one if
$\mtx{A}$ is exactly rank-one deficient, as in our case.
Since $\vct{y}$ is in the range of~$\mtx{A}$, we have that
$\mtx{A}\vct{x} = \vct{y}$, and therefore $\vct{z} = \vct{x} -
\vct{r}$ is a null-vector of $\mtx{A}$.
We then normalize $\vct{z}$, which represents a discretization of
$\sigma$, the generalized Debye source, so that $\int_\gamma
|\sigma|^2 ds = 1$. Components of the zero-flux Beltrami fields for
the first five $\lambda$'s greater than give are show in
Figure~\ref{fig_eigs}.

\section{Conclusions}
\label{sec_conclusion}

We have presented a numerical algorithm for the calculation of
Beltrami (force-free) fields in the interior of toroidally
axisymmetric geometries, both genus one and two, and given explicit
numerical examples demonstrating the accuracy of the solver. The
representation of Beltrami fields is based on earlier work regarding
the representation of Maxwell fields in the context of scattering
phenomena. It ensures that (by construction) the resulting field is a
Beltrami field, and leads directly to a second-kind integral equation
which can be numerically inverted to high precision. Moreover, the
solver has low memory requirements since the unknowns in the integral
equation are defined on the boundary of the domain; the interior of
the domain does not need to be discretized. The solver can be used to
study toroidally axisymmetric Taylor states in laboratory, space, and
astrophysical plasmas, as well as bifurcations to the fully
three-dimensional, non-axisymmetric states corresponding to resonances
of the Beltrami parameter $\lambda$.

Our formulation extends directly to non-axisymmetric domains, as
required for stellarator equilibrium calculations for example
\cite{hudson}. In that context, the results presented here are
promising, since high accuracy is reached with a modest number of
unknowns, even for low aspect ratio, highly elongated domains. This,
along with the low memory requirements, are desirable features for
computations of three-dimensional MHD equilibria which rely on the
iterative calculation of multiple Taylor states in genus-one and
genus-two regions \cite{hudson}.  The fact that the interior of the
domain is not discretized in our scheme and that none of the
quantities need to be evaluated in the interior is particularly
attractive for that application, since the force balance conditions
between adjacent Taylor states only need to be evaluated at the ideal
MHD interfaces at each iteration. The entire iterative procedure could
thus be implemented by discretizing the ideal MHD surfaces only, and
the magnetic field would only be evaluated in the entire domain at the
very end, once global force balance has been reached. This would lead
to significant savings in run-time and in memory, and motivates
detailed code-to-code comparisons with the scheme used in SPEC \cite{hudson} to identify the strengths and weaknesses of each approach. However, for our numerical solver to be applicable to
general three dimensional toroidal domains, details in the numerical
implementation need to be modified. Specifically, we are currently
working on addressing questions regarding the discretization of the
domain, the evaluation of differential surface operators, and the
computation of numerical quadratures.

Furthermore, regarding code-to-code comparisons, one thing seems clear
at this stage: the scheme of this work is the only one in which
only the boundary is discretized. This leads to a major reduction in
the number of unknowns as compared to other existing schemes, as we
have highlighted on several occasions throughout the paper.

Lastly, with minor changes the representation and algorithm described
are applicable in exterior, unbounded geometries (both simply- and
multiply-connected). We have, however, left this discussion out of
this article because of the present focus on plasma physics
applications.

\bibliographystyle{abbrv}
\bibliography{preprint}

\end{document}